\documentstyle[]{article}

\newfont{\bb}{msbm10}

\def\Bbb#1{\hbox{\bb#1}}

\def\func#1{\hbox{\tt#1}}
\def\overleftrightarrow#1{\stackrel{\leftrightarrow}{#1}}

\def\text#1{#1}
\newtheorem{theorem}{Theorem}[section]

\newtheorem{corollary}{Corollary}[section]
\newtheorem{conjecture}{Conjecture}[section]

\begin{document}

\title{Induced surfaces and their integrable dynamics. II. 
Generalized Weierstrass
representations in 4D spaces and deformations via DS hierarchy.}
\author{B. G. Konopelchenko and G. Landolfi\\
\em{Dipartimento di Fisica, Universit\`a di Lecce, 73100 Lecce, Italy} }
\date{}
\maketitle

\begin{abstract}
Extensions of the generalized Weierstrass representation to generic surfaces
in 4D Euclidean and pseudo-Euclidean spaces are given. Geometric
characteristics of surfaces are calculated. It is shown that integrable
deformations of such induced surfaces are generated by the Davey-Stewartson
hierarchy. Geometrically these deformations are characterized by the
invariance of an infinite set of functionals over surface. The Willmore
functional (the total squared mean curvature) is the simplest of them.
Various particular classes of surfaces and their integrable deformations are
considered.
\end{abstract}

\newpage

\section{Introduction}
\setcounter{equation}{0}

Surfaces and their deformations (dynamics) were for a long time the 
subjects of
intensive study both in mathematics and physics. Theory of immersions and
deformations of surfaces has been a significant part of the classical
differential geometry (see \textit{e.g.} [1-3]). Various methods to describe
immersions and deformations have been developed. This subject continues to
be an important part of the contemporary differential geometry (see 
\textit{e.g.} [4-6]).

In physics, the dynamics of interfaces, surfaces, fronts is a key ingredient
in a number of interesting phenomena from hydrodynamics, propagation of
flame fronts, growth of crystals, deformations of membranes to world-sheets
and their dynamics in the string theory (see \textit{e.g.} [7-9]). Such
dynamics could be modelled by nonlinear partial differential equations.
Analytic methods to study surfaces, their properties and deformations are of
great interest both in mathematics and physics (see [7-9] and recent papers
[10-11]).

A general method to construct surfaces via the solutions of linear
differential equations and their deformations via the corresponding
nonlinear integrable equations has been proposed in [12-13]. The two basic
examples considered in [13] were given by: \emph{1)} the generalized
Weierstrass formulae for generic surfaces conformally immersed into $\Bbb{R}%
^{3}$ and deformations via the modified Veselov-Novikov equation; \emph{2)}
the Lelieuvre formula for surfaces in $\Bbb{R}^{3}$ referred to asymptotic
lines and integrable deformations via the Nizhnik-Veselov-Novikov equation.

The generalized Weierstrass representation proposed in [12-13] has been
proved to be an effective tool to study generic surfaces in $\Bbb{R}^{3}$
and their deformations. In differential geometry its use has allowed to
obtain several interesting results both of local and global character, in
particular, for the Willmore functional $W=\int H^{2}\left[ dS\right] $
where $H$ is the mean curvature of the surface (see \emph{e.g.} [17-22]). In
physics, it has been applied to study of various problems in the theory of
liquid membranes, 2D gravity and string theory [17, 23-27]. In the string
theory the functional $W=\int H^{2}\left[ dS\right] $ is known as the
Polyakov extrinsic action and in membrane theory it is the Helfrich free
energy [7-9]. Deformations constructed via the Lelieuvre formula and the
Nizhnik-Veselov equation have occurred to be an interesting class of
deformations of surfaces in affine geometry [28].

An extension of the Weierstrass representation to multidimensional spaces
would be of a great interest. In physics, a strong motivation lies in the
Polyakov string integral over surfaces in multidimensional spaces [7-9].
Theory of immersion of surfaces into four-dimensional spaces is an important
part of the contemporary differential geometry too [4-6, 29-32].

In this paper we present extensions of the generalized Weierstrass
representation to the cases of generic surfaces conformally immersed into
the four-dimensional spaces $\Bbb{R}^{4}$, $\Bbb{R}^{3,1}$ and $\Bbb{R}%
^{2,2} $ with the metrics $g_{ik}=diag(1,1,1,1)$, $g_{ik}=diag(1,1,1,-1)$
and $g_{ik}=diag(1,1,-1,-1)$, respectively. A basic linear system consists
of a couple of the two-dimensional Dirac equations while the formulae for
immersions, the induced metric, mean curvature and Willmore functional are
of the type similar to that of the $\Bbb{R}^{3}$ case.

Integrable deformations of surfaces are generated by the Davey-Stewartson
(DS) hierarchy of 2+1-dimensional soliton equations. These deformations of
surfaces inherit all remarkable properties of the soliton equations.
Geometrically, such deformations are characterized by the invariance of an
infinite set of functionals over surfaces. The simplest of them is given by
the Willmore functional.

We consider both space-like and time-like surfaces in the case of
pseudo-Euclidean spaces. Various particular classes of surfaces, including
minimal and superminimal surfaces, immersions with constant mean curvature
and their deformations are discussed. It is shown that one special class of
the 1+1-dimensional reductions gives rise to an integrable motions of curves
on the three-dimensional sphere $\Bbb{S}^{3}$ and hyperboloids described by
the nonlinear Schroedinger equation.

Extensions of the generalized Weierstrass formulae to the four-dimensional
Riemann spaces are discussed too.

The present paper can be considered as the second part of the paper [13]. In
fact, the DS inducing of surfaces (in $\Bbb{R}^{3}$) and their deformations
via the DS hierarchy have been mentioned in [13] (section 10, pp. 41-42) as
one of possibilities. Here we elaborate this case in detail.

Note that some results of this paper has been presented briefly in [33]. The
Weierstrass type representations for particular classes of surfaces have
been discussed in [30, 34] and recently in [35].

The paper is organized as follows. In section 2 we remind for convenience
some results concerning the generalized Weierstrass formulae for surfaces in 
$\Bbb{R}^{3}$. In section 3 we present the Weierstrass representation for
generic surfaces conformally immersed into $\Bbb{R}^{4}$. Space-like
surfaces in $\Bbb{R}^{2,2}$, $\Bbb{R}^{3,1}$ are considered in section 4
while in section 5 we concentrate on time-like surfaces in pseudo-Euclidean
spaces. Surfaces in four-dimensional Riemann spaces are discussed in section
6. Integrable deformations of surfaces in 4D spaces via the DS hierarchy are
described in section 7. Explicit expressions for deformations of surfaces
are given in section 8. Particular classes of surfaces and their
deformations are considered in section 9. The one-dimensional reduction of
the Weierstrass representation and corresponding surfaces are discussed in
section 10. Appendix contains some basic facts about the DS hierarchy.

\section{Generalized Weierstrass formulae for surfaces in $\Bbb{R}^{3}$.}
\setcounter{equation}{0}

A generalization of the Weierstrass formulae to generic surfaces in $\Bbb{R}%
^{3}$ proposed by one of the authors in 1993 (see [12] and [13]) starts with
the linear system (two-dimensional Dirac equation) 
\begin{equation}
\begin{array}{l}
\psi _{z}=p\varphi \quad , \\ 
\varphi _{\overline{z}}=-p\psi
\end{array}
\label{2.1}
\end{equation}
where $\psi $ and $\varphi $ are complex-valued functions of $z$, $\overline{%
z}\in \Bbb{C}$ and $p\left( z,\overline{z}\right) $ is a real-valued
function. Then one defines the three real-valued functions $X^{1}\left( z,%
\overline{z}\right) $, $X^{2}\left( z,\overline{z}\right) $ and $X^{3}\left(
z,\overline{z}\right) $ by the formulae 
\begin{eqnarray}
X^{1}+iX^{2}=i\int_{\Gamma }(\overline{\psi }^{2}dz^{\prime }-\overline{%
\varphi }^{2}d\overline{z}^{\prime })\quad , \nonumber \\ 
X^{1}-iX^{2}=i\int_{\Gamma }(\varphi ^{2}dz^{\prime }-\psi ^{2}d\overline{z}%
^{\prime })\quad , \nonumber \\ 
X^{3}=-\int_{\Gamma }(\overline{\psi }\varphi \,dz^{\prime }+\psi \overline{%
\varphi \,}d\overline{z}^{\prime })\quad 
\label{2.2}
\end{eqnarray}
where $\Gamma $ is an arbitrary contour in $\Bbb{C}$. In virtue of (\ref{2.1}%
) the \emph{r.h.s.} in (\ref{2.2}) do not depend on the choice of $\Gamma $.
If one now treats $X^{i}\left( z,\overline{z}\right) $ as the coordinates in 
$\Bbb{R}^{3}$ then the formulae (\ref{2.1}), (\ref{2.2}) define a conformal
immersion of surface into $\Bbb{R}^{3}$ with the induced metric of the form 
\begin{equation}
ds^{2}=u^{2}\,dz\,d\overline{z}=\left( |\psi |^{2}+|\varphi |^{2}\right)
^{2}dz\,d\overline{z}\qquad ,  \label{2.3}
\end{equation}
with the Gauss curvature 
\begin{equation}
K=-\frac{4}{u^{2}}\left[ \log u\right] _{z\overline{z}}  \label{2.4}
\end{equation}
and the mean curvature 
\begin{equation}
H=2\frac{p}{u}  \quad . \label{2.5} 
\end{equation}
At $p=0$ one gets minimal surfaces and the formulae (\ref{2.2}) are reduced
to the old Weierstrass formulae.

Another analog of the Weierstrass formulae for surfaces of prescribed (non
zero) mean curvature have been proposed earlier by Kenmotsu in [36]. The
Kenmotsu representation is given by 
\begin{equation}
\overrightarrow{X}=\func{Re}\left[ \int_{\Gamma }\eta \overrightarrow{\phi }%
\,dz^{\prime }\right]  \label{2.7}
\end{equation}
where $\overrightarrow{\phi }=\left[ 1-f^{2},i\left( 1+f^{2}\right)
,2f\right] $ and the functions $f$ and $\eta $ obey the following
compatibility condition 
\begin{equation}
\left( \log \eta \right) _{\overline{z}}=-\frac{2\overline{f}f_{\overline{z}}%
}{1+\left| f\right| ^{2}}\quad .  \label{2.8}
\end{equation}
Here and below the bar denotes the complex conjugation. Then the mean
curvature $H$ is 
\begin{equation}
H=-\frac{2f_{\overline{z}}}{\overline{\eta }\left( 1+\left| f\right|
^{2}\right) ^{2}}\quad .  \label{2.9}
\end{equation}
It $\,$was proved in [9] that any surface in $\Bbb{R}^{3}$ can be presented
in such a form. This representation of surfaces deals basically with the
Gauss map for generic surface in $\Bbb{R}^{3}$ [31].

It turned out that the Kenmotsu formulae (\ref{2.7}), (\ref{2.8}) and the
generalized Weierstrass formulae (\ref{2.1}), (\ref{2.2}) are equivalent to
each other. The relation between the functions $(f$,$\eta )$ and $(\psi
,\varphi )$ is the following [14] 
\begin{equation}
f=i\frac{\overline{\psi }}{\varphi }\qquad ,\qquad \eta =i\varphi ^{2}
\label{2.10}
\end{equation}
and 
\begin{equation}
p=-\frac{\eta f_{\overline{z}}}{\sqrt{\eta \overline{\eta }}\left( 1+\left|
f\right| ^{2}\right) }  \quad . \label{2.11}
\end{equation}
So all results proved for the Kenmotsu formulae [36] and associated Gauss
map [31] in $\Bbb{R}^{3}$ are valid also for the generalized Weierstrass
formulae (\ref{2.2}). In particular, it implies immediately that any surface 
$\Bbb{R}^{3}$ in can be represented via (\ref{2.1})-(\ref{2.2}).

Though the representations (\ref{2.1}), (\ref{2.2}) and (\ref{2.7}), (\ref
{2.8}) are equivalent, the former provides us certain advantages. They are
mainly due to the fact that in the generalized Weierstrass formulae the
functions $\psi $ and $\varphi $ obey the 
linear equations (\ref{2.1}) while for
the Kenmotsu formulae the nonlinear constraint (\ref{2.8}) is 
difficult to deal with.
This circumstance has allowed to simplify essentially an analysis that had
lead to several interesting results both of local and global character
[15-21]. It occurred, in particular, that the Willmore functional (see \emph{%
e.g.} [4]) or the Helfrich-Polyakov action (see [7-9]) $W=\int 
\overrightarrow{H}^{2}\left[ dS\right] $ has a very simple form: $W=4\int
p^{2}dx\,dy$ ($z=x+iy$) [14-15].

One of the advantages of the generalized Weierstrass formulae (\ref{2.1}), (%
\ref{2.2}) is that they allow to construct a new class of deformations of
surfaces via the modified Veselov-Novikov equation [12-13]. The
characteristic feature of these integrable deformations is that the Willmore
functional remains invariant [14-15]. Thus, the generalized Weierstrass
representation (\ref{2.2}) has been proved to be an effective tool to study
surfaces in $\Bbb{R}^{3}$ and their deformations.

We would like to 
emphasize that the idea to generate surfaces via solutions of
linear equations is, in fact, the old idea of the classical differential
geometry as it was already noted in [13]. In [3] one can find the two
representations of these type in addition to the Weierstrass formulae. The
first is given by the Lelieuvre formula

\begin{equation}
\begin{array}{l}
\overrightarrow{X_{\xi }}=\overrightarrow{\nu }\times \overrightarrow{\nu
_{\xi }}\quad , \\ 
\overrightarrow{X_{\eta }}=-\overrightarrow{\nu }\times \overrightarrow{\nu
_{\eta }}
\end{array}
\label{2.12}
\end{equation}
where $\overrightarrow{\nu }=\left( \nu _{1},\nu _{2},\nu _{3}\right) $ are
three linearly independent solutions of the equation

\begin{equation}
\overrightarrow{\nu _{\xi \eta }}+p\overrightarrow{\nu }=0  \label{2.13}
\end{equation}
and the $\nu _{i}(\xi ,\eta )$ 's and $p(\xi ,\eta )$ are scalar functions.
The Lelieuvre formula defines immersion of a surface into $\Bbb{R}^{3}$ ($%
\overrightarrow{X}:\Bbb{R}^{2}\rightarrow \Bbb{R}^{3}$) parametrized by
asymptotic lines $\xi =const$ and $\eta =const$ . The Lelieuvre formula is
well-known in affine geometry of surfaces.

Another example is provided by the equation 
\begin{equation}
\theta _{\eta \xi }-(\log \lambda )_{\eta }\theta _{\xi }-\lambda ^{2}\theta
=0  \label{2.14}
\end{equation}
where $\xi $ and $\eta $ are real variables and $\lambda $ is a real-valued
function. It is stated in [3] that two solutions of (\ref{2.14}) define, via
certain integral formulae, a surface in $\Bbb{R}^{3}$ parametrized by
minimal lines, but no calculation of the metric and curvature is given. This
example, seems, was forgotten completely until it has been found during the
preparation of the second paper [13] on the generalized Weierstrass
formulae. The representation (\ref{2.14}) is rather close to that of (\ref
{2.1})-(\ref{2.2}). Indeed, equation (\ref{2.14}) can be rewritten as the
system 
\begin{equation}
\begin{array}{l}
\theta _{\xi }=\lambda \varphi \quad , \\ 
\varphi _{\eta }=\lambda \theta
\end{array}
\label{2.15}
\end{equation}
where $\varphi $ is a new function. If one takes two solutions $(\theta
,\varphi )$ and $(\widetilde{\theta },\widetilde{\varphi })$ of the system (%
\ref{2.15}) then the formulae given in [3] (pp. 82) take the form 
\begin{eqnarray}
X^{1}+iX^{2}=\int (\theta ^{2}d\eta +\varphi ^{2}d\xi )\quad , 
\nonumber \\ 
X^{1}-iX^{2}=\int (\widetilde{\theta }^{2}d\eta +\widetilde{\varphi }%
^{2}d\xi )\quad ,  \nonumber \\ 
X^{3}=i\int (\theta \widetilde{\theta }\,d\eta +\varphi \widetilde{\varphi }%
\,d\xi )\quad .
\label{2.16}
\end{eqnarray}
However, in contrast to the representation (\ref{2.1}), (\ref{2.2}), the
formulae (\ref{2.15}), (\ref{2.16}) do not define a real surface in $\Bbb{R}%
^{3}$.

We would like to note that some results in [37] and [38] were close to the
generalized Weierstrass representation (\ref{2.2}). In [37] a formula
similar to (\ref{2.2}) for constant mean curvature surfaces has been
discussed. In [38] the system (\ref{2.1}) had appeared within the
quaternionic description of surfaces in $\Bbb{R}^{3}$ (formula (2.19) of
[38]). However in [38] it was accompanied by another two equations (equation
(2.16) of [38]) which are indispensable in the Sym's type approach. So the
meaning of the system (\ref{2.1}), seems, has been missed. The generalized
Weierstrass type formulae admit also a beautiful formulation within the
spinor representation of surfaces [39-40].

\section{The Weierstrass representation for surfaces in $\Bbb{R}^{4}$}
\setcounter{equation}{0}

An extension of the representation (\ref{2.1}), (\ref{2.2}) to the
four-dimensional Euclidean space is quite natural. 
Let $\psi _{1}$, $\varphi
_{1}$ and $\psi _{2}$, $\varphi _{2}$ be solutions of the systems 
\begin{equation}
\begin{array}{l}
\psi _{1z}=p\varphi _{1}\quad , \\ 
\varphi _{1\overline{z}}=-\overline{p}\psi _{1}
\end{array}
\qquad 
\begin{array}{l}
\psi _{2z}=\overline{p}\varphi _{2}\quad , \\ 
\varphi _{2\overline{z}}=-p\psi _{2}\qquad .
\end{array}
\label{3.1}
\end{equation}
Equations (\ref{3.1}) imply that 
\begin{equation}
\left( \psi _{1}\psi _{2}\right) _{z}=-\left( \varphi _{1}\varphi
_{2}\right) _{\overline{z}}\qquad ,\qquad \left( \psi _{1}\overline{\varphi }%
_{2}\right) _{z}=\left( \varphi _{1}\overline{\psi }_{2}\right) _{\overline{z%
}}\quad .  \label{3.2}
\end{equation}
As a consequence there are four functions $X^{i}\left( z,\overline{z}\right) 
$ ($i=1,2,3,4$) such that 
\begin{eqnarray}
dX^{1}=\frac{1}{2}\left( \overline{\psi }_{1}\overline{\psi }_{2}-\varphi
_{1}\varphi _{2}\right) dz+c.c. \nonumber \\ 
dX^{2}=\frac{i}{2}\left( \overline{\psi }_{1}\overline{\psi }_{2}+\varphi
_{1}\varphi _{2}\right) dz+c.c. \nonumber \\ 
dX^{3}=\frac{1}{2}\left( \varphi _{1}\overline{\psi }_{2}+\psi _{1}\overline{%
\varphi }_{2}\right) dz+c.c. \nonumber \\ 
dX^{4}=\frac{i}{2}\left( \overline{\psi }_{1}\varphi _{2}-\varphi _{1}%
\overline{\psi }_{2}\right) dz+c.c.
\label{3.3}
\end{eqnarray}
where $c.c.$ denotes a complex conjugation of the previous term. We treat
now these functions $X^{i}(z,\overline{z})$ as the coordinates of surfaces
in $\Bbb{R}^{4}$. For components of induced metric 
\begin{equation}
g_{zz}=\sum_{i=1}^{4}\left( X_{z}^{i}\right) ^{2}=\overline{g_{\overline{z}\,%
\overline{z}}}\qquad ,\qquad g_{z\overline{z}}=\sum_{i=1}^{4}\left(
X_{z}^{i}X_{\overline{z}}^{i}\right)  \label{3.4}
\end{equation}
one gets 
\begin{equation}
g_{zz}=g_{\overline{z}\,\overline{z}}=0  \label{3.5}
\end{equation}
and 
\begin{equation}
g_{z\overline{z}}=\frac{1}{2}\left( \left| \psi _{1}\right| ^{2}+\left|
\varphi _{1}\right| ^{2}\right) \left( \left| \psi _{2}\right| ^{2}+\left|
\varphi _{2}\right| ^{2}\right) \quad .  \label{3.6}
\end{equation}
Further, two normal vectors $\overrightarrow{N}_{1}$, $\overrightarrow{N}%
_{2} $ are 
\begin{equation}
\overrightarrow{N}_{1}=\sqrt{\frac{\left| \varphi _{1}\right| ^{2}\left|
\varphi _{2}\right| ^{2}}{u_{1}u_{2}}}\func{Re}\left( \overrightarrow{A}%
\right) \quad ,\quad \overrightarrow{N}_{2}=\sqrt{\frac{\left| \varphi
_{1}\right| ^{2}\left| \varphi _{2}\right| ^{2}}{u_{1}u_{2}}}\func{Im}\left( 
\overrightarrow{A}\right)  \label{3.7}
\end{equation}
where 
\begin{equation}
\begin{array}{l}
u_{k}=\left| \psi _{k}\right| ^{2}+\left| \varphi _{k}\right| ^{2}\qquad
,\qquad k=1,2 \\ \vspace{-5mm} \\
\overrightarrow{A}=\left[ -\frac{\psi _{1}}{\overline{\varphi }_{1}}-\frac{%
\overline{\psi }_{2}}{\varphi _{2}},i\left( \frac{\psi _{1}}{\overline{%
\varphi }_{1}}-\frac{\overline{\psi }_{2}}{\varphi _{2}}\right) ,\frac{\psi
_{1}\overline{\psi }_{2}}{\overline{\varphi }_{1}\varphi _{2}}-1,-i\left( 1+%
\frac{\psi _{1}\overline{\psi }_{2}}{\overline{\varphi }_{1}\varphi _{2}}%
\right) \right]
\end{array}
\quad .  \label{3.8-3.9}
\end{equation}
The mean curvature vector defined standardly as 
\[
\overrightarrow{H}=\frac{1}{g_{z\overline{\,z}}}\overrightarrow{X}_{z%
\overline{\,z}} 
\]
is given by 
\begin{eqnarray}
\vec{H} &=&\frac{2}{u_{1}u_{2}}\left[ Re\left( p\varphi _{1} \psi_2+%
\overline{p}\psi _{1}\varphi _{2}\right) ,Im\left( p\varphi _{1}
\psi_2+%
\overline{p}\psi _{1}\varphi _{2}\right) ,\right.  \nonumber \\
&&\left. \;\;\;\;\;\;\;\;\;\;\;\;\;Re\left( p\varphi _{1}\overline{\varphi }%
_{2}-\overline{p}\psi _{1}\overline{\psi }_{2}\right) ,Im\left( p\overline{%
\psi }_{1}\psi _{2}-\overline{p} \overline{\varphi }_{1}\varphi _{2}\right)
\right] \; .  \label{3.10}
\end{eqnarray}
The components $h_{1}$, $h_{2}$ of $\overrightarrow{H}$ along $%
\overrightarrow{N}_{1}$ and $\overrightarrow{N}_{2}$ 
$\left( \text{\emph{i.e.%
} }\overrightarrow{H}=h_{1}\overrightarrow{N}_{1}+h_{2}\overrightarrow{N}%
_{2}\right) $are 
\begin{equation}
h_{1}=-\frac{2\func{Re}\left( p\overline{\varphi }_{1}\varphi _{2}\right) }{%
\sqrt{\,|\varphi _{1}|^{2}|\varphi _{2}|^{2}\,u_{1}u_{2}}}\quad \quad
,\qquad h_{2}=\frac{2\func{Im}\left( p\overline{\varphi }_{1}\varphi
_{2}\right) }{\sqrt{\,|\varphi _{1}|^{2}|\varphi _{2}|^{2}\,u_{1}u_{2}}}
\quad . \label{3.11}
\end{equation}
So, the mean curvature $\overrightarrow{H}^{2}=%
\sum_{i=1}^{4}H^{i}H^{i}=h_{1}^{2}+h_{2}^{2}$ is equal to 
\begin{equation}
\overrightarrow{H}^{2}=4\frac{\left| p\right| ^{2}}{u_{1}u_{2}}\quad .
\label{3.12}
\end{equation}
Then the Gaussian curvature is 
\begin{equation}
K=-\frac{2}{u_{1}u_{2}}\left[ \log \left( u_{1}u_{2}\right) \right] _{z%
\overline{z}}\quad .  \label{3.13}
\end{equation}
Finally, the Willmore functional $W=\int \vec{H}^{2}\left[ dS\right] $ is
given by 
\begin{equation}
W=4\int \left| p\right| ^{2}dx\,dy\text{\thinspace \quad }.  \label{3.14}
\end{equation}
Thus, we have the following

\begin{theorem}
The generalized Weierstrass formulae 
\begin{eqnarray}
X^{1}+iX^{2} &=&\int_{\Gamma }{\left( -\varphi _{1}\varphi _{2}dz^{\prime
}+\psi _{1}\psi _{2}d\overline{z}^{\prime }\right) }\;,  \nonumber \\
X^{1}-iX^{2} &=&\int_{\Gamma }{\left( \overline{\psi }_{1}\overline{\psi }%
_{2}dz^{\prime }-\overline{\varphi }_{1}\overline{\varphi }_{2}d\overline{z}%
^{\prime }\right) }\;,  \nonumber \\
X^{3}+iX^{4} &=&\int_{\Gamma }{\left( \varphi _{1}\overline{\psi }%
_{2}dz^{\prime }+\psi _{1}\overline{\varphi }_{2}d\overline{z}^{\prime
}\right) }\;,  \nonumber \\
X^{3}-iX^{4} &=&\int_{\Gamma }{\left( \overline{\psi }_{1}\varphi
_{2}dz^{\prime }+\overline{\varphi }_{1}\psi _{2}d\overline{z}^{\prime
}\right) }  \label{3.15}
\end{eqnarray}
where 
\begin{equation}
\begin{array}{l}
\psi _{1z}=p\varphi _{1}\quad , \\ 
\varphi _{1\overline{z}}=-\overline{p}\psi _{1}\quad
\end{array}
\quad \quad ,\qquad 
\begin{array}{l}
\psi _{2z}=\overline{p}\varphi _{2}\quad , \\ 
\varphi _{2\overline{z}}=-p\psi _{2} \quad,
\end{array}
\label{3.16}
\end{equation}
$\Gamma $ is a contour in $\Bbb{C}$, define the conformal immersion of a
surface into $\Bbb{R}^{4}$. The induced metric is of the form 
\begin{equation}
ds^{2}=u_{1}u_{2}\,dz\,d\overline{z}  \label{3.17}
\end{equation}
where $u_{k}=\left| \psi _{k}\right| ^{2}+\left| \varphi _{k}\right| ^{2}$ ($%
k=1,2$), the Gaussian and squared mean curvatures are 
\begin{equation}
K=-\frac{2}{u_{1}u_{2}}\left[ \log \left( u_{1}u_{2}\right) \right] _{z%
\overline{z}}\quad \quad ,\quad \quad \overrightarrow{H}^{2}=4\frac{\left|
p\right| ^{2}}{u_{1}u_{2}}\quad .  \label{3.18}
\end{equation}
The total squared mean curvature (Willmore functional) is given by 
\begin{equation}
W=4\int \left| p\right| ^{2}dx\,dy\text{\thinspace \quad .}  \label{3.19}
\end{equation}
\end{theorem}

Since the linear system (\ref{3.16}) contains two arbitrary functions ($%
\func{Re}p$ and $\func{Im}p$) of two variables, then the formulae (\ref{3.15}%
), (\ref{3.16}) allows us to get any surface in $\Bbb{R}^{4}$. The
generalized Weierstrass representation (\ref{3.15}) defines surfaces in $%
\Bbb{R}^{4}$ up to translations. In the specialized case $\overline{p}=p$
one gets the formulae derived in [35]. In the particular case $\psi _{2}=\pm
\psi _{1}$, $\varphi _{2}=\pm \varphi _{1}$, $X_{z}^{4}=X_{\overline{z}%
}^{4}=0$ and the formulae (\ref{3.16})-(\ref{3.19}) are reduced to those (%
\ref{2.1}), (\ref{2.2}) of the $\Bbb{R}^{3}$ case with the substitution 
$X^{1}\leftrightarrow X^{2}$, $X^{3}\leftrightarrow -X^{3}$.

Note that a linear system of the form (\ref{3.1}) arises also as the
restriction of the Dirac equation to a surface in $\Bbb{R}^{4}$ [41].

Note the equations (\ref{3.3}) can be represented in the form 
\begin{eqnarray}
d\left( X^{1}+iX^{2}\right)& = &
-\varphi _{1}\varphi _{2}dz+\psi _{1}\psi _{2}d%
\overline{z}\quad , \nonumber \\ 
d\left( X^{3}+iX^{4}\right) & = &
\varphi _{1}\overline{\psi }_{2}dz+\psi _{1}%
\overline{\varphi }_{2}d\overline{z}
\label{3.20}
\end{eqnarray}
which reveals a symmetry between the pairs of coordinates $(X^{1},X^{2})$
and $(X^{3},X^{4})$.

The formulae (\ref{3.3}) can be also rewritten in a spinor representation
type form 
\[
d\left( \sigma _{1}X^{1}+\sigma _{2}X^{2}+\sigma _{3}X^{3}+iIX^{4}\right)
=V_{2}^{\dagger }\left( 
\begin{array}{ll}
0 & dz \\ 
d\overline{z} & 0
\end{array}
\right) V_{1} 
\]
where 
\[
V_{1,2}=\left( 
\begin{array}{cc}
\psi _{1,2} & -\overline{\varphi }_{1,2} \\ 
\varphi _{1,2} & \overline{\psi }_{1,2}
\end{array}
\right) \quad , 
\]
$\sigma _{i}$ ($i=1,2,3$) are the Pauli matrices and $I$ is the identity
matrix.

The condition (\ref{3.5}), that an immersion is conformal, written as 
\begin{equation}
(X_{z}^{1})^{2}+(X_{z}^{2})^{2}+(X_{z}^{2})^{2}+(X_{z}^{2})^{2}=0
\label{3.24}
\end{equation}
defines the complex quadric $\Bbb{Q}_{2}$%
\begin{equation}
w_{1}^{2}+w_{2}^{2}+w_{3}^{2}+w_{4}^{2}=0  \label{3.25}
\end{equation}
in $\Bbb{CP}^{3}$ where $w_{i}$ ($i=1,2,3,4$) are homogeneous coordinates. A
diffeomorphism of $\Bbb{Q}_{2}$ to the Grassmannian $\Bbb{G}_{2,4}$ of
oriented 2-planes in $\Bbb{R}^{4}$ allows us to define the Gauss map $%
\overrightarrow{G}(z)$ for a surface represented by the generalized
Weierstrass formulae (\ref{3.15}). It is given by 
\begin{equation}
\overrightarrow{G}(z)=\left[ 1+\frac{\overline{\psi }_{1}\overline{\psi }_{2}%
}{\varphi _{1}\varphi _{2}},i\left( 1-\frac{\overline{\psi }_{1}\overline{%
\psi }_{2}}{\varphi _{1}\varphi _{2}}\right) ,i\left( \frac{\overline{\psi }%
_{1}}{\varphi _{1}}+\frac{\overline{\psi }_{2}}{\varphi _{2}}\right) ,\left( 
\frac{\overline{\psi }_{1}}{\varphi _{1}}-\frac{\overline{\psi }_{2}}{%
\varphi _{2}}\right) \right]  \quad . \label{3.26}
\end{equation}
The Gauss map for surfaces immersed into $\Bbb{R}^{4}$has been studied
earlier in the paper [31]. In [31] the Gauss map $\overrightarrow{G}(z)$ has
been parametrized as follows 
\begin{equation}
\overrightarrow{G}(z)=\left[ 1+f_{1}f_{2},i\left( 1-f_{1}f_{2}\right)
,f_{1}-f_{2},-i\left( f_{1}+f_{2}\right) \right]  \label{3.27}
\end{equation}
where $f_{1}$ and $f_{2}$ are complex-valued functions. A surface in $\Bbb{R}%
^{4}$ is then defined by [31] 
\begin{equation}
\overrightarrow{X}=\int_{\Gamma }\func{Re}\left( \eta \overrightarrow{G}%
\,dz\right)  \label{3.28}
\end{equation}
where $f_{1}$ and $f_{2}$ satisfy the compatibility conditions 
\begin{equation}
\func{Im}\left[ \left( \frac{f_{1z\overline{z}}}{f_{1\overline{z}}}-2\frac{%
\overline{f_{1}}f_{1z}}{1+\left| f_{1}\right| ^{2}}\right) _{\overline{z}%
}+\left( \frac{f_{2z\overline{z}}}{f_{2\overline{z}}}-2\frac{\overline{f_{2}}%
f_{2}}{1+\left| f_{2}\right| ^{2}}\right) _{\overline{z}}\right] =0
\label{3.29}
\end{equation}
and 
\begin{equation}
\left| F_{1}\right| =\left| F_{2}\right|  \label{3.30}
\end{equation}
where $F_{i}=f_{i\overline{z}}\left( 1+\left| f_{i}\right| ^{2}\right) ^{-1}$%
, $i=1,2$. The function $\eta $ is given by 
\begin{equation}
\overline{\eta }^{2}=-\frac{4F_{1}F_{2}}{H^{2}\left( 1+\left| f_{1}\right|
^{2}\right) \left( 1+\left| f_{2}\right| ^{2}\right) }  \label{3.31}
\end{equation}
where the mean curvature $H$ is expressed via $f_{1}$ and $f_{2}$ by 
\begin{equation}
2\left( \log H\right) _{z}=\frac{f_{1z\overline{z}}}{f_{1\overline{z}}}-2%
\frac{\overline{f_{1}}f_{1z}}{1+\left| f_{1}\right| ^{2}}+\frac{f_{2z%
\overline{z}}}{f_{2\overline{z}}}-2\frac{\overline{f_{2}}f_{2}}{1+\left|
f_{2}\right| ^{2}}\quad .  \label{3.32}
\end{equation}
Similar to the three-dimensional case this representation includes the
complicated compatibility conditions.

\begin{theorem}
The generalized Weierstrass representation (\ref{3.15})-(\ref{3.19}) and the
Gauss map type representation (\ref{3.28})-(\ref{3.32}) are equivalent to
each other via the substitution 
\begin{equation}
\eta =i\varphi _{1}\varphi _{2}\quad ,\quad f_{1}=i\frac{\overline{\psi }_{1}%
}{\varphi _{1}}\quad ,\quad f_{2}=-i\frac{\overline{\psi }_{2}}{\varphi _{2}}%
\quad .
\end{equation}
\label{333}
\end{theorem}

The proof is straightforward: equations (\ref{3.16}) and (\ref{333}) give
the constraints (\ref{3.29})-(\ref{3.31}) with 
\begin{equation}
p=-iF_{1}\frac{\varphi _{1}}{\overline{\varphi }_{1}}\qquad ,\qquad 
\overline{p}=iF_{2}\frac{\varphi _{2}}{\overline{\varphi }_{2}}  \label{3.34}
\end{equation}
while (\ref{3.15}) is converted into (\ref{3.28}).

\section{Generalized Weierstrass representations for surfaces in
pseudo-Euclidean spaces}
\setcounter{equation}{0}

The derivation of the generalized Weierstrass formulae for surfaces immersed
into four-dimensional pseudo-Euclidean spaces is rather similar to that of $%
\Bbb{R}^{4}$.

\begin{theorem}
The generalized Weierstrass formulae 
\begin{eqnarray}
X^{1}+iX^{2} &=&\int_{\Gamma }(\varphi _{1}\varphi _{2}\,dz^{\prime }+\psi
_{1}\psi _{2}d\overline{z}^{\prime })\;,  \nonumber \\
X^{1}-iX^{2} &=&\int_{\Gamma }(\,\overline{\psi }_{1}\overline{\psi }%
_{2}dz^{\prime }+\overline{\varphi }_{1}\overline{\varphi }_{2}d\overline{z}%
^{\prime })\;,  \nonumber \\
X^{3}+iX^{4} &=&i\int_{\Gamma }(\,\overline{\psi }_{1}\varphi _{2}dz^{\prime
}+\overline{\varphi }_{1}\psi _{2}d\overline{z}^{\prime })\;,  \nonumber \\
X^{3}-iX^{4} &=&-i\int_{\Gamma }(\varphi _{1}\,\overline{\psi }%
_{2}dz^{\prime }+\psi _{1}\overline{\varphi }_{2}d\overline{z}^{\prime })
\label{4.1}
\end{eqnarray}
where 
\begin{equation}
\begin{array}{l}
\psi _{1z}=p\varphi _{1}\;\;, \\ 
\varphi _{1\overline{z}}=\overline{p}\psi _{1}
\end{array}
\begin{array}{l}
\psi _{2z}=\overline{p}\varphi _{2}\;, \\ 
\varphi _{2\overline{z}}=p\psi _{2}\;\;\;,
\end{array}
\label{4.2}
\end{equation}
$\psi _{\alpha }$, $\varphi _{\alpha }$, $p$ are complex-valued functions, $%
\Gamma $ is a contour in $\Bbb{C}$, define the conformal immersion $\vec{X}:%
\Bbb{C}\rightarrow \Bbb{R}^{2,2}$ of a surface into the space $\Bbb{R}^{2,2}$%
. The induced metric is 
\begin{equation}
ds^{2}=v_{1}v_{2}dzd\overline{z}  \label{4.3}
\end{equation}
where $v_{\alpha }=|\psi _{\alpha }|^{2}-|\varphi _{\alpha }|^{2}$, $\alpha
=1,2$, the Gaussian and mean curvature are of the form 
\begin{equation}
K=-\frac{2}{v_{1}v_{2}}\left[ \log {\left( v_{1}v_{2}\right) }\right] _{z%
\overline{z}}\;\;,\;\;\vec{H}^{2}=-\frac{4|p|^{2}}{v_{1}v_{2}}\;\;
\label{4.4}
\end{equation}
and the Willmore functional $W=\int {\vec{H}^{2}\left[ dS\right] }$ is given
by 
\begin{equation}
W=-4\int \left| p\right| ^{2}dx\,dy\;\;.  \label{4.5}
\end{equation}
\end{theorem}

The proof is similar to the case of $\Bbb{R}^{4}$, only now equations (%
\ref{4.2}) give $(\psi _{1}\psi _{2})_{z}=(\varphi _{1}\varphi _{2})_{%
\overline{z}}$ and $(\psi _{1}\overline{\varphi }_{2})_{z}=(\varphi _{1}%
\overline{\psi }_{2})_{\overline{z}}$. The mean curvature vector is given by 
\begin{eqnarray}
\vec{H} &=&\frac{2}{v_{1}v_{2}}\left[ \func{Re}\left( p\varphi _{1} 
\psi_2+%
\overline{p}\psi _{1}\varphi _{2}\right) ,
\func{Im}\left( p\varphi_{1} \psi_2
+\overline{p}\psi _{1}\varphi _{2}\right) ,\right.  \nonumber \\
&&\left. \;\;\;\;\;\;\;\;\;\;\;\;\;\func{Im}\left( p\varphi _{1}\overline{%
\varphi }_{2}+\overline{p}\psi _{1}\overline{\psi }_{2}\right) ,\func{Re}%
\left( p\overline{\psi }_{1}\psi _{2}+\overline{p}\overline{\varphi }%
_{1}\varphi _{2}\right) \right] \;.  \label{4.6}
\end{eqnarray}
In the particular case $\overline{p}=p$ the formulae (\ref{4.1})-(\ref{4.5})
are reduced to those obtained in [{35}] with the substitution $%
X^{1}\leftrightarrow X^{2}$, $X^{3}\leftrightarrow X^{4}$. In contrast to
[35] the formulae (\ref{4.1})-(\ref{4.5}) allow to represent an arbitrary
surface in $\Bbb{R}^{2,2}$.

Conformal immersions into the Minkowski space $\Bbb{R}^{3,1}$ are given by
slightly different formulae.

\begin{theorem}
The Weierstrass type formulae 
\begin{eqnarray}
X^{1}+iX^{2}=\int_{\Gamma }(\varphi _{1}\overline{\psi }_{2}dz^{\prime
}+\psi _{1}\overline{\varphi }_{2}d\overline{z}^{\prime })\quad , 
\nonumber \\ 
X^{1}-iX^{2}=\int_{\Gamma }(\overline{\psi }_{1}\varphi _{2}dz^{\prime }+%
\overline{\varphi }_{1}\psi _{2}d\overline{z}^{\prime })\quad , 
\nonumber \\ 
X^{3}+X^{4}=\int_{\Gamma }(\overline{\psi }_{1}\varphi _{1}dz^{\prime }+\psi
_{1}\overline{\varphi }_{1}d\overline{z}^{\prime })\quad , 
\nonumber \\ 
X^{3}-X^{4}=-\int_{\Gamma }(\overline{\psi }_{2}\varphi _{2}dz^{\prime
}+\psi _{2}\overline{\varphi }_{2}d\overline{z}^{\prime })
\label{4.7}
\end{eqnarray}
where 
\begin{equation}
\begin{array}{l}
\psi _{\alpha z}=p\varphi _{\alpha }\;, \\ 
\varphi _{\alpha \overline{z}}=q\psi _{\alpha }
\end{array}
\quad ,\quad \alpha =1,2  \label{4.8}
\end{equation}
$q$ and $p$ are real-valued functions, $\Gamma $ is a contour in $\Bbb{C}$,
define the conformal immersion of a surface into the Minkowski space $%
\overrightarrow{X}:\Bbb{C}\rightarrow \Bbb{R}^{3,1}$. 
The induced metric on a
surface is 
\begin{equation}
ds^{2}=\left| \psi _{1}\varphi _{2}-
\varphi _{1} \psi_2 \right| ^{2}dz\,d%
\overline{z}\quad ,  \label{4.9}
\end{equation}
the Gaussian curvature is 
\[
K=-\frac{2}{|\psi _{1}\varphi _{2}-\psi _{2}\varphi _{1}|^{2}}\,\left[ \log
\left( |\psi _{1}\varphi _{2}-\varphi _{1} \psi_2|\right) \right] _{z\,%
\overline{z}} 
\]
the squared mean curvature $\overrightarrow{H}^{2}$ and the Willmore
functional are given respectively by 
\begin{equation}
\overrightarrow{H}^{2}=-\frac{4qp}{\left| \psi _{1}\varphi _{2}-\psi
_{2}\varphi _{1}\right| ^{2}}\quad ,\quad W=-4\int qp\,dx\,dy\quad .
\label{4.10}
\end{equation}
\end{theorem}

In this case the linear system (\ref{4.8}) implies that 
\[
\left( \psi _{\alpha }\overline{\varphi }_{\beta }\right) _{z}=\left(
\varphi _{\alpha }\overline{\psi }_{\beta }\right) _{\overline{z}%
}\;\;\;\;\alpha ,\beta =1,2 
\]
that guarantee an independence of the \emph{r.h.s.} of (\ref{4.7}) on the
choice of the contour $\Gamma $ of integration. The rest is straightforward.
In particular, the mean curvature vector $\overrightarrow{H}$ is of the form 
\begin{eqnarray}
\vec{H} &=&\frac{2}{\left| \psi _{1}\varphi _{2}-\varphi_{1} \psi_2
\right| ^{2}}\left[ \func{Re}\left( p\varphi _{1}\overline{\varphi }%
_{2}+q\psi _{1}\overline{\psi }_{2}\right) ,\func{Re}\left( ip\overline{%
\varphi }_{1}\varphi _{2}+iq\overline{\psi }_{1}\psi _{2}\right) ,\right. 
\nonumber \\
&&
\left. \frac{p}{2}\left( \left| \varphi
_{1}\right| ^{2}-\left| \varphi _{2}\right| ^{2}\right) +\;\frac{q}{2}\left(
\left| \psi _{1}\right| ^{2}-\left| \psi _{2}\right| ^{2}\right),  
\frac{p}{2}\left( \left|
\varphi _{1}\right| ^{2}+\left| \varphi _{2}\right| ^{2}\right) +\;\frac{q}{2%
}\left( \left| \psi _{1}\right| ^{2}+\left| \psi _{2}\right| ^{2}\right)
\;\right] \;. \nonumber \\
&&
\label{4.11}
\end{eqnarray}
Since again one has two arbitrary real-valued functions $p$ and $q$, the
Weierstrass type formulae (\ref{4.7})-(\ref{4.8}) allow us to construct any
surface immersed into $\Bbb{R}^{3,1}$.

Differential version of all three generalized Weierstrass representations
given above can be written in the following common form 
\begin{equation}
d\left( \sum_{i=1}^{4}{\tau _{i}X^{i}}\right) =\Phi _{2}^{\dagger }\left( 
\begin{array}{cc}
0 & dz \\ 
d\overline{z} & 0
\end{array}
\right) \Phi _{1}  \label{4.12}
\end{equation}
where $\dagger $ denotes Hermitian conjugation. In the case of immersion
into $\Bbb{R}^{4}$ one has 
\[
\tau _{1}=\sigma _{1}\;\;,\;\;\tau _{2}=\sigma _{2}\;\;,\;\;\tau _{3}=\sigma
_{3}\;\;,\;\;\tau _{4}=i\sigma _{4} 
\]
and 
\begin{equation}
\Phi _{\alpha }=\left( 
\begin{array}{cc}
\psi _{\alpha } & -\overline{\varphi }_{\alpha } \\ 
\varphi _{\alpha } & \overline{\psi }_{\alpha }
\end{array}
\right) \;\;,\;\;\alpha =1,2  \label{4.13}
\end{equation}
where $\sigma _{1}$, $\sigma _{2}$, $\sigma _{3}$ are the standard Pauli
matrices and $\sigma _{4}$ is an identical $2\times 2$ matrix. At the $\Bbb{R%
}^{2,2}$ case 
\[
\tau _{1}=\sigma _{1}\;\;,\;\;\tau _{2}=\sigma _{2}\;\;,\;\;\tau
_{3}=i\sigma _{3}\;\;,\;\;\tau _{4}=\sigma _{4} 
\]
and 
\begin{equation}
\Phi _{\alpha }=\left( 
\begin{array}{cc}
\psi _{\alpha } & \overline{\varphi }_{\alpha } \\ 
\varphi _{\alpha } & \overline{\psi }_{\alpha }
\end{array}
\right) \;\;,\;\;\alpha =1,2  \quad .\label{4.14}
\end{equation}
Finally, the immersion into the Minkowski space $\Bbb{R}^{3,1}$ correspond
to 
\[
\tau _{i}=\sigma _{i}\;\;\;\;(i=1,2,3,4) 
\]
and 
\begin{equation}
\Phi _{1}=\Phi _{2}=\left( 
\begin{array}{cc}
\psi _{1} & \psi _{2} \\ 
\varphi _{2} & \varphi _{2}
\end{array}
\right) \;\;.  \label{4.15}
\end{equation}
In fact, one can start with the formulae (\ref{4.12}) to derive the
Weierstrass representations in the forms (\ref{3.15})-(\ref{3.16}), (\ref
{4.1})-(\ref{4.2}) and (\ref{4.7})-(\ref{4.8}). Indeed, one can show that
the $1-$form in the \emph{r.h.s.} of (\ref{4.12}) is closed if the $2\times
2 $ matrices $\Phi _{1}$, $\Phi _{2}$ obey the Dirac equations 
\begin{equation}
\left( 
\begin{array}{cc}
\partial _{z} & 0 \\ 
0 & \partial _{\overline{z}}
\end{array}
\right) \Phi _{1}=\left( 
\begin{array}{cc}
u & p \\ 
q & v
\end{array}
\right) \Phi _{1}\;\;,\;\;\left( 
\begin{array}{cc}
\partial _{z} & 0 \\ 
0 & \partial _{\overline{z}}
\end{array}
\right) \Phi _{2}=\left( 
\begin{array}{cc}
-\overline{v} & \overline{p} \\ 
\overline{q} & -\overline{u}
\end{array}
\right) \Phi _{2}  \label{4.16}
\end{equation}
where $p$, $q$, $u$, $v$ are arbitrary complex-valued functions. Functions $%
u $ and $v$ always can be converted to zeros by gauge transformation
(redefinition of $\Phi $). Then the reality conditions for $X^{i}$ are
satisfied if matrices $\Phi _{\alpha }$ have the form (\ref{4.13}), (\ref
{4.14}) or (\ref{4.15}) while the functions $p$, $q$ should obey the
constraints $p+\overline{q}=0$, $p-\overline{q}=0$ and $p=\overline{p}$, $q=%
\overline{q}$, respectively. Consequently, the corresponding formula (\ref
{4.12}) gives rise to the Weierstrass representations considered above.

A formula of the type (\ref{4.12}) appears naturally [42] in the
quaternionic approach to surfaces (see also [{38-39,43}]) which could
provide an invariant formulation of the construction presented above.

In the particular case $\overline{p}=p$, $\psi _{1}=\psi _{2}=\psi $ and $%
\varphi _{1}=\varphi _{2}=\varphi $ the Weierstrass representation (\ref{4.1}%
), (\ref{4.2}) defines the conformal immersion 
\begin{equation}
\begin{array}{l}
X^{1}+iX^{2}=\int_{\Gamma }(\varphi ^{2}dz^{\prime }+\psi ^{2}d\overline{z}%
^{\prime }) \\ \vspace{-5mm} \\
X^{4}=\int_{\Gamma }(\overline{\psi }\varphi dz^{\prime }+\psi \overline{%
\varphi }d\overline{z}^{\prime })
\end{array}
\label{4.17}
\end{equation}
where 
\begin{equation}
\begin{array}{l}
\psi _{z}=p\varphi \;, \\ 
\varphi _{\overline{z}}=p\psi
\end{array}
\label{4.18}
\end{equation}
of a surface into the three-dimensional pseudo-Euclidean space with the
metric $g_{ik}=diag(1,1,-1)$. The induced metric is 
\begin{equation}
ds^{2}=\left( \left| \psi \right| ^{2}-\left| \varphi \right| ^{2}\right)
^{2}\,dz\,d\overline{z}  \label{4.19}
\end{equation}
while the squared mean curvature and the Willmore functional are 
\begin{equation}
\overrightarrow{H}^{2}=-\frac{4p^{2}}{\left( \left| \psi \right| ^{2}-\left|
\varphi \right| ^{2}\right) ^{2}}\qquad ,\qquad W=-4\int p^{2}dx\,dy
\label{4.20}
\end{equation}
respectively.

So (\ref{4.17}), (\ref{4.18}) give the pseudo-Euclidean version of the
generalized representation (\ref{2.1}), (\ref{2.2}). This representation
again generates any surface conformally immersed in $\Bbb{R}^{2,1}$. The
representation (\ref{4.17}), (\ref{4.18}) for surfaces in $\Bbb{R}^{2,1}$
has been found earlier in [44]. The Kenmotsu type formula for nonminimal
surfaces in $\Bbb{R}^{2,1}$ has been given in [45]. An analog of the
formulae (\ref{2.7})-(\ref{2.9}) is of the form [45] 
\begin{equation}
\overrightarrow{X}=\func{Re}\left( \int^{z}\eta \overrightarrow{\phi }%
dz^{\prime }\right)  \label{4.21}
\end{equation}
where 
\begin{equation}
\overrightarrow{\phi }=\left[ 1+f^{2},i\left( 1-f^{2}\right) ,2f\right]
\label{4.22}
\end{equation}
and 
\begin{equation}
2\overline{f}_{z}=\eta H\left( 1-|f|^{2}\right) ^{2}\quad ,\quad (\log
H)_{z}=\frac{1}{f_{\overline{z}}}\left( f_{z\overline{z}}+2\frac{\overline{f}%
f_{z}f_{\overline{z}}}{1-|f|^{2}}\right) \quad . \label{4.23}
\end{equation}
Comparing the Gauss map vector for the Weierstrass type representation (\ref
{4.17}), (\ref{4.18}), \emph{i.e.} 
\begin{equation}
G\left( z\right) =\left[ \left( 1+\frac{\varphi ^{2}}{\overline{\psi }^{2}}%
\right) ,i\left( 1-\frac{\varphi ^{2}}{\overline{\psi }^{2}}\right) ,2\frac{%
\varphi }{\overline{\psi }}\right]  \label{4.25}
\end{equation}
with (\ref{4.21}), (\ref{4.22}), one concludes that 
\begin{equation}
f=\frac{\varphi }{\overline{\psi }}\qquad ,\qquad \eta =\overline{\psi }^{2}
\quad . \label{4.26}
\end{equation}
Further, the relation 
\[
H=\frac{2p}{\left| \psi \right| ^{2}-\left| \varphi \right| ^{2}} 
\]
converts the nonlinear equations (\ref{4.23}) into the linear system (\ref{4.18}%
) and vice versa.

So the Weierstrass representation (\ref{4.17}), (\ref{4.18}) and the
Kenmotsu type formulae (\ref{4.21})-(\ref{4.23}) are equivalent to each
other through the relations (\ref{4.26}). Obviously, the Kenmotsu type
formulae (\ref{4.21})-(\ref{4.23}) can be obtained also as a particular case
of the 
generalized Weierstrass representation for surfaces immersed into $\Bbb{R}%
^{3,1}$. Indeed, the formulae (\ref{4.17})-(\ref{4.18}) arise
from the 
formulae (\ref{4.7})-(\ref{4.8}) under the reduction $p=q,\overline{%
\psi }_{2}=\varphi _{1},\varphi _{2}=\overline{\psi }_{1}$.

\section{The Weierstrass representations for time-like surfaces in the 4D
pseudo-Euclidean spaces}
\setcounter{equation}{0}

Hyperbolic (time-like) surfaces with the signature $(+,-)$ appear naturally
in the pseudo-Euclidean spaces. They arise as the world-sheets for strings
moving in Minkowski space or other pseudo-Euclidean spaces (see [7-9]) and
are of interest in differential geometry too [5,6]. For time-like surfaces
minimal lines are real. So, to get a surface parametrized by minimal lines,
one has to start with the Dirac linear equations which instead of $z$, $%
\overline{z}$ contain real independent variables (say $\xi $ and $\eta $).

\begin{theorem}
The formulae 
\begin{equation}
\begin{array}{l}
X^{1}=\frac{1}{2}\int_{\Gamma }\left[ \left( \overline{\varphi }_{1}\varphi
_{2}+\varphi _{1}\overline{\varphi }_{2}\right) d\xi' +\left( \overline{\psi }%
_{1}\psi _{2}+\psi _{1}\overline{\psi }_{2}\right) d\eta' \right] \quad , 
\\ \vspace{-5mm} \\ 
X^{2}=\frac{i}{2}\int_{\Gamma }\left[ \left( \overline{\varphi }_{1}\varphi
_{2}-\varphi _{1}\overline{\varphi }_{2}\right) d\xi' +\left( \overline{\psi }%
_{1}\psi _{2}-\psi _{1}\overline{\psi }_{2}\right) d\eta' \right] \quad , 
\\ \vspace{-5mm} \\ 
X^{3}=\frac{1}{2}\int_{\Gamma }\left[ \left( 
\varphi _{1} \overline{\varphi}_1-
\varphi_{2} \overline{\varphi}_2 \right) d\xi' +
\left( \psi _{1} \overline{\psi}_1 -\psi _{2} \overline{\psi}_2 \right) 
d\eta' \right] \quad , \\ \vspace{-5mm} \\
X^{4}=\frac{1}{2}\int_{\Gamma }\left[ \left( 
\varphi _{1} \overline{\varphi}_1+
\varphi _{2} \overline{\varphi}_2 \right) d\xi' +
\left( 
\psi _{1} \overline{\psi}_1+\psi _{2} \overline{\psi}_2
\right) d\eta' \right] \quad ,
\end{array}
\label{5.1}
\end{equation}
where the functions $\psi _{\alpha }$, $\varphi _{\alpha }$ ($\alpha =1,2$)
satisfy the linear system 
\begin{equation}
\begin{array}{l}
\psi _{\alpha \xi }=p\varphi _{\alpha }\quad , \\ 
\varphi _{\alpha \eta }=\overline{p}\psi _{\alpha }\quad ,
\end{array}
\qquad \alpha =1,2  \label{5.2}
\end{equation}
$\Gamma $ is a contour of integration in $\Bbb{R}^{2}$, define an immersion
into the Minkowski space $\Bbb{R}^{3,1}$of a generic time-like surface
parametrized by the minimal lines $\xi =const$, $\eta =const$ .The induced
metric is 
\begin{equation}
ds^{2}=-\left| \varphi _{1}\psi _{2}-\psi _{1}\varphi _{2}\right| ^{2}d\eta
\,d\xi \quad ,  \label{5.3}
\end{equation}
the Gaussian and squared mean curvature are given respectively by 
\begin{equation}
K=\frac{2}{|\varphi _{1}\psi _{2}-\psi _{1}\varphi _{2}|^{2}}\left[ \log
\left( |\varphi _{1}\psi _{2}-\psi _{1}\varphi _{2}|\right) \right] _{z\,%
\overline{z}}\quad ,\quad \overrightarrow{H}^{2}=4\frac{|p|^{2}}{|\varphi
_{1}\psi _{2}-\psi _{1}\varphi _{2}|^{2}}  \label{5.4}
\end{equation}
while the Willmore functional is 
\begin{equation}
W=2\int |p|^{2}d\eta \,d\xi  \quad . \label{5.5}
\end{equation}
\end{theorem}

The formulae (\ref{5.1}), (\ref{5.2}) allow to construct any surface in $%
\Bbb{R}^{3,1}$. The mean curvature vector defined as 
\[
\overrightarrow{H}=\frac{1}{g_{\eta \xi }}\overrightarrow{X}_{\eta \xi } 
\]
is of the form 
\begin{eqnarray}
\overrightarrow{H} &=&-\frac{|\varphi _{1}\psi _{2}-\psi _{1}\varphi
_{2}|^{-2}}{2}\func{Re}[p(\overline{\psi }_{1}\varphi _{2}+\overline{\psi }%
_{2}\varphi _{1}),ip(\overline{\psi }_{1}\varphi _{2}+\overline{\psi }%
_{2}\varphi _{1}),  \nonumber \\
&&\qquad \qquad \qquad \qquad \qquad \qquad p(\overline{\psi }_{1}\varphi
_{1}-\overline{\psi }_{2}\varphi _{2}),p(\overline{\psi }_{1}\varphi _{1}+%
\overline{\psi }_{2}\varphi _{2})] \;.  \nonumber \\
&&
\label{5.6}
\end{eqnarray}
It is interesting that the formulae similar to (\ref{5.1}) have appeared in
completely different context as parametrization of constraint for
coordinates of classical string governed by the Nambu-Goto Lagrangian [46].
The formulae from [46] look like 
\begin{equation}
\begin{array}{l}
X_{\xi }^{1}=ac+bd\quad , \\ 
X_{\xi }^{2}=ad-bc\quad , \\ 
\left( X^{4}+X^{3}\right) _{\xi }=a^{2}+b^{2}\quad , \\ 
\left( X^{4}-X^{3}\right) _{\xi }=c^{2}+d^{2}
\end{array}
\label{5.8}
\end{equation}
where $a$, $b$, $c$, $d$ are real-valued function; similar formulae with
real-valued functions $\widetilde{a},\widetilde{b},\widetilde{c},\widetilde{%
d}$ hold for $X_{\eta }^{i}$. It is not difficult to see that (\ref{5.8})
and (\ref{5.1}) coincide under the identification 
\begin{equation}
\begin{array}{lll}
\varphi _{1}=a+ib\quad , & \qquad & \varphi _{2}=c+id \\ 
\psi _{1}=\widetilde{a}+i\widetilde{b}\quad , & \qquad & \psi _{2}=%
\widetilde{c}+i\widetilde{d}
\end{array}
\label{5.9}
\end{equation}
Note that the paper [46] contains no the linear system (\ref{5.2}) and
formulae (\ref{5.3})-(\ref{5.6}). The result of the paper [46] suggest
possible applications of the representation (\ref{5.1}), (\ref{5.2}) in the
theory of classical strings. For example, the Nambu-Goto action takes the
form 
\[
S_{NG}=\frac{\alpha }{2}\int \left| \varphi _{1}\psi _{2}-\psi _{1}\varphi
_{2}\right| ^{2}d\eta \,d\xi \quad . 
\]

\begin{theorem}
The Weierstrass type formulae 
\begin{equation}
\begin{array}{l}
X^{1}=\frac{1}{2}\int_{\Gamma }[(\varphi _{1}\widetilde{\varphi }_{2}+%
\widetilde{\varphi }_{1}\varphi _{2})d\xi ^{\prime }+(\psi _{1}\widetilde{%
\psi }_{2}+\widetilde{\psi }_{1}\psi _{2})d\eta ^{\prime }]\quad , 
\\ \vspace{-5mm} \\ 
X^{2}=\frac{1}{2}\int_{\Gamma }[(\varphi _{1}\widetilde{\varphi }_{1}-%
\widetilde{\varphi }_{2}\varphi _{2})d\xi ^{\prime }+(\psi _{1}\widetilde{%
\psi }_{1}-\psi _{2}\widetilde{\psi }_{2})d\eta ^{\prime }]\quad , 
\\ \vspace{-5mm} \\ 
X^{3}=\frac{1}{2}\int_{\Gamma }[(\widetilde{\varphi }_{1}\varphi
_{2}-\varphi _{1}\widetilde{\varphi }_{2})d\xi ^{\prime }+(\widetilde{\psi }%
_{1}\psi _{2}-\psi _{1}\widetilde{\psi }_{2})d\eta ^{\prime }]\quad , 
\\ \vspace{-5mm} \\ 
X^{4}=\frac{1}{2}\int_{\Gamma }[(\varphi _{1}\widetilde{\varphi }%
_{1}+\varphi _{2}\widetilde{\varphi }_{2})d\xi ^{\prime }+(\psi _{1}%
\widetilde{\psi }_{1}+\widetilde{\psi }_{2}\psi _{2})d\eta ^{\prime }]
\end{array}
\label{5.10}
\end{equation}
where 
\begin{equation}
\begin{array}{lll}
\psi _{\alpha \xi }=p\varphi _{\alpha }\quad , & \quad \quad & \widetilde{%
\psi }_{\alpha \xi }=q\widetilde{\varphi }_{\alpha }\quad , \\ 
\varphi _{\alpha \eta }=q\psi _{\alpha }\quad , & \quad \quad & \widetilde{%
\varphi }_{\alpha \eta }=p\widetilde{\psi }_{\alpha }\quad ,
\end{array}
\qquad (\alpha =1,2)  \label{5.11}
\end{equation}
$\psi _{\alpha }$, $\varphi _{\alpha }$, $\widetilde{\psi }_{\alpha }$, $%
\widetilde{\varphi }_{\alpha }$ $(\alpha =1,2)$ and $p$, $q$ are real-valued
functions, $\Gamma $ is a contour of integration in $\Bbb{R}^{2}$, define an
immersion into $\Bbb{R}^{2,2}$ of a generic surface parametrized by minimal
lines. The induced metric is 
\begin{equation}
ds^{2}=(\psi _{1}\varphi _{2}-\psi _{2}\varphi _{1})(\widetilde{\varphi }_{1}%
\widetilde{\psi }_{2}-\widetilde{\psi }_{1}\widetilde{\varphi }_{2})d\xi
\,d\eta  \label{5.02}
\end{equation}
while the Gaussian and the squared mean curvatures are 
\begin{equation}
K=-\frac{2}{(\psi _{1}\varphi _{2}-\psi _{2}\varphi _{1})(\widetilde{\varphi 
}_{1}\widetilde{\psi }_{2}-\widetilde{\psi }_{1}\widetilde{\varphi }_{2})}%
\left\{ \log \left[ (\psi _{1}\varphi _{2}-\psi _{2}\varphi _{1})(\widetilde{%
\varphi }_{1}\widetilde{\psi }_{2}-\widetilde{\psi }_{1}\widetilde{\varphi }%
_{2})\right] \right\} _{z\,\overline{z}} 
\end{equation}
\begin{equation}
\overrightarrow{H}^{2}=\frac{4qp}{(\psi _{1}\varphi _{2}-\psi _{2}\varphi
_{1})(\widetilde{\psi }_{1}\widetilde{\varphi }_{2}-\widetilde{\psi }_{2}%
\widetilde{\varphi }_{1})}
\end{equation}
and the Willmore functional is given by 
\begin{equation}
W=-2\epsilon \int qp\,d\xi \,d\eta \quad .  \label{5.14}
\end{equation}
where $\epsilon =sign\left[ (\psi _{1}\varphi _{2}-\psi _{2}\varphi _{1})(%
\widetilde{\varphi }_{1}\widetilde{\psi }_{2}-\widetilde{\psi }_{1}%
\widetilde{\varphi }_{2})\right] $.
\end{theorem}

In this case one has from (\ref{5.11}) 
\[
\left( \psi _{\alpha }\widetilde{\psi }_{\beta }\right) _{\xi }=\left(
\varphi _{\alpha }\widetilde{\varphi }_{\beta }\right) _{\eta }\qquad
,\qquad \alpha ,\beta =1,2 
\]
that gives rise to (\ref{5.10}). The mean curvature vector is 
\[
\begin{array}{ll}
\overrightarrow{H} & =\frac{1}{g_{\eta \xi }}\overrightarrow{X_{\eta \xi }}%
=\left[ (\psi _{1}\varphi _{2}-\psi _{2}\varphi _{1})(\widetilde{\varphi }%
_{1}\widetilde{\psi }_{2}-\widetilde{\psi }_{1}\widetilde{\varphi }%
_{2})\right] ^{-1}[(q\psi _{1}\widetilde{\varphi }_{2}+p\varphi _{1}%
\widetilde{\psi }_{2})+(1\leftrightarrow 2), \\ 
& (q\psi _{1}\widetilde{\varphi }_{1}+p\varphi _{1}\widetilde{\psi }%
_{1})-(1\leftrightarrow 2),(q\psi _{2}\widetilde{\varphi }_{1}+p\psi _{2}%
\widetilde{\varphi }_{1})-(1\leftrightarrow 2),(q\psi _{1}\widetilde{\varphi 
}_{1}+p\varphi _{1}\widetilde{\psi }_{1})+(1\leftrightarrow 2)].
\end{array}
\]

For time-like surfaces an analog of the formula (\ref{4.12})\ looks like 
\[
d\left( \sum \tau _{i}X^{i}\right) =\Phi _{2}^{\dagger }\left( 
\begin{array}{ll}
0 & d\xi \\ 
d\eta & 0
\end{array}
\right) \Phi _{1} 
\]
where for the case of the space $\Bbb{R}^{3,1}${} one has 
\[
\tau _{1}=\sigma _{1}\quad ,\quad \tau _{2}=i\sigma _{2}\quad ,\quad \tau
_{3}=\sigma _{3}\quad ,\quad \tau _{4}=\sigma _{4} 
\]
and 
\[
\Phi _{1}=\Phi _{2}=\left( 
\begin{array}{cc}
\psi _{1} & \psi _{2} \\ 
\varphi _{1} & \varphi _{2}
\end{array}
\right) 
\]
while for the space $\Bbb{R}^{2,2}$%
\[
\tau _{1}=\sigma _{1}\quad ,\quad \tau _{2}=\sigma _{3}\quad ,\quad \tau
_{3}=i\sigma _{2}\quad ,\quad \tau _{4}=\sigma _{4} 
\]
and 
\[
\Phi _{1}=\left( 
\begin{array}{ll}
\psi _{1} & \psi _{2} \\ 
\varphi _{1} & \varphi _{2}
\end{array}
\right) \qquad ,\qquad \Phi _{2}=\left( 
\begin{array}{ll}
\widetilde{\varphi }_{1} & \widetilde{\varphi }_{2} \\ 
\widetilde{\psi }_{1} & \widetilde{\psi }_{2}
\end{array}
\right) \quad . 
\]
Note that in the special case $\psi _{1}=\psi _{2}$, $\varphi _{1}=\varphi
_{2}$ (for (\ref{5.1})-(\ref{5.2})) and at $\psi _{1}=\psi _{2}$, $\varphi
_{1}=\varphi _{2}$, $\widetilde{\psi }_{1}=\widetilde{\psi }_{2}$, $%
\widetilde{\varphi }_{1}=\widetilde{\varphi }_{2}$ (for (\ref{5.10})-(\ref
{5.11}) ) the immersions (\ref{5.1})-(\ref{5.2}) and (\ref{5.10})-(\ref{5.11}%
) are reduced to maps $\Bbb{R}^{2}\rightarrow \Bbb{R}^{2}$ ( $(\xi ,\eta
)\rightarrow (X^{1},X^{4})$ or $(\xi ,\eta )\rightarrow (X^{1},X^{2})$ ).

It is interesting that the generalized Weierstrass formulae given in the
theorems (3.1), (4.1), (4.2), (5.1), (5.2) 
take place also in the case when the quantities $\psi _{\alpha} $ 
and $\varphi _{\alpha }$ are elements of the Grassmannian algebra, 
\textit{i.e.} when they anticommute to each other. 
The geometric characteristics of surfaces are given by formulae similar to
those of the theorems (3.1), (4.1), (4.2), (5.1), (5.2). 
This type of the Weierstrass representation is of
the interest in the string theory.

\section{Surfaces in four-dimensional Riemann space}
\setcounter{equation}{0}

The results of the previous sections apparently can be extended to the case
of immersion into generic four-dimensional Riemann space with the metric
tensor $g_{ik}$.

\begin{theorem}
The generalized Weierstrass formulae define an immersion of surface into the
four dimensional Riemann space with the metric tensor $g_{ik}$. The induced
metric is 
\begin{equation}
ds^{2}=g_{\xi \xi }\,d\xi ^{2}+2g_{\xi \eta }\,d\xi \,d\eta +g_{\eta \eta
}\,d\eta ^{2}  \label{6.1}
\end{equation}
where 
\begin{equation}
g_{\xi \xi }=g_{ik}X_{\xi }^{i}X_{\xi }^{k}\;\;,\;\;g_{\xi \eta
}=g_{ik}X_{\xi }^{i}X_{\eta }^{k}\quad ,\quad g_{\eta \eta }=g_{ik}X_{\eta
}^{i}X_{\eta }^{k}\quad .  \label{6.2}
\end{equation}
\end{theorem}

The Gaussian and mean curvature are calculated straightforwardly. One can
choose any of the formulae presented above to define the coordinates $X^{i}$
via the solutions of the linear system. It is obvious, however, that to get
an immersion which is, for instance, locally conformal around a point one
should choose the Weierstrass type formula for the 4D pseudo-Euclidean space
with the metric which coincide with the signature of the desired 4D Riemann
space.

In the case of conformally-Euclidean spaces 
($g_{ik}=e^{2\sigma }\delta _{ik}$, 
$i,k=1,2,3,4$, $\sigma $ is a function and $\delta _{ik}$ is the Kronecker
symbol) the immersion is the conformal one: 
\begin{equation}
ds^{2}=e^{2\sigma }u_{1}u_{2}dzd\overline{z}\;\;.  \label{6.3}
\end{equation}
The Gaussian and mean curvatures are 
\begin{equation}
K=-2e^{-2\sigma }\frac{\left[ 2\sigma +\log {\left( u_{1}u_{2}\right) }%
\right] _{z\overline{z}}}{u_{1}u_{2}}\;\;,\;\;H^{2}=4e^{-2\sigma }\frac{%
|p|^{2}}{u_{1}u_{2}}\;\;.  \label{6.4}
\end{equation}
For the Willmore functional one gets 
\begin{equation}
W=4\int |p|^{2}dx\,dy\;\;.  \label{6.5}
\end{equation}
A special case of immersions into the space $\Bbb{S}^{4}$ of constant
curvature had attracted recently the particular interest (see \emph{e.g.} [{4%
},{47}]). To describe it we choose the Riemann form for the metric of $\Bbb{S%
}^{4}$, \emph{i.e.} (see \emph{e.g.} [{48}]) 
\begin{equation}
e^{2\sigma }=\left[ 1+\frac{K_{0}}{4}\sum_{i=1}^{4}\left( X^{i}\right)
^{2}\right] ^{-2}  \label{6.6}
\end{equation}
where $K_{0}$ is the curvature. Then the formulae (\ref{3.15}), (\ref{3.16})
define the conformal immersion of a surface into $\Bbb{S}^{4}$. At $\psi
_{2}=\pm \psi _{1}$, $\varphi _{2}=\pm \varphi _{1}$ ($X^{4}=const$) 
one has the
conformal immersion into $\Bbb{S}^{3}$. The generalized Weierstrass
representation provides us an effective method to study immersions into $%
\Bbb{S}^{3}$ and $\Bbb{S}^{4}$ and, consequently, the Willmore surfaces 
in $\Bbb{R}^3$.

Similarly, one can get conformal immersions into hyperbolic spaces $\Bbb{S}%
^{3,1}$ and $\Bbb{S}^{2,2}$ with the constant curvature. Indeed, one has
simply to take 
\begin{equation}
e^{2\sigma }=\left[ 1+\frac{K_{0}}{4}\overrightarrow{X}^{2}\right]
^{-2}\quad .  \label{6.7}
\end{equation}

\section{Integrable deformations}
\setcounter{equation}{0}

In construction of deformations of surfaces given by the Weierstrass
representations we follow to the general approach of [12-13].

So we assume that all quantities in the linear problems (\ref{3.16}), (\ref
{4.2}), (\ref{4.8}), (\ref{5.2}), (\ref{5.11}) (except $z$, 
$\overline{z}$, $\xi$ and $\eta $)
depend on the new deformations parameters $t_{n}$. Then we assume that this
dependence on $t_{n}$ is such that there are operators $A_{n}$, $B_{n}$, $%
C_{n}$, $D_{n}$ such that equations 
\begin{equation}
\begin{array}{l}
\psi _{t_{n}}=A_{n}\psi +B_{n}\varphi \quad , \\ 
\varphi _{t_{n}}=A_{n}\psi +B_{n}\varphi
\end{array}
\qquad ,\qquad n=1,2,3,\ldots  \label{7.1}
\end{equation}
hold. The compatibility conditions of (\ref{3.16}), (\ref{4.2}), (\ref{4.8}%
), (\ref{5.2}), (\ref{5.11}) with (\ref{7.1}) fix the dependence of $\psi $, 
$\varphi $ and $p$, $q$ on $t_{n}$ and, consequently, define the
deformations of surfaces. Concrete cases are governed by different
specializations (reductions) of the DS hierarchy.

Let us consider first immersions of space-like surfaces into the Minkowski
space $\Bbb{R}^{3,1}$. In this case $p$ and $q$ are real-valued functions.
The corresponding deformations are generated by the ''real'' DSII hierarchy
(with real- valued $p$ and $q$). In particular, in equations (A.2), (A.5)
the constants $\alpha $, $\gamma $, $\alpha _{2}$ are real. Since 
(\ref{A.12}) is
obviously the integral of motion also for real-valued $p$ and $q$, then, in
virtue of (\ref{4.10}), the Willmore functional $W$ remains invariant under
these DSII deformations.

For the Weierstrass representations of the space-like surfaces in $\Bbb{R}%
^{4}$ and $\Bbb{R}^{2,2}$ we have the linear system (A.1) with the following
reductions: 
\begin{equation}
\left( 
\begin{array}{cc}
0 & p_{1} \\ 
q_{1} & 0
\end{array}
\right) =\left( 
\begin{array}{cc}
0 & p \\ 
\varepsilon \overline{p} & 0
\end{array}
\right) \;\;\;\;for\;\;\;\Phi _{1}  \label{7.2}
\end{equation}
and 
\begin{equation}
\left( 
\begin{array}{cc}
0 & p_{2} \\ 
q_{2} & 0
\end{array}
\right) =\left( 
\begin{array}{cc}
0 & \overline{p} \\ 
\varepsilon p & 0
\end{array}
\right) \;\;\;\;for\;\;\;\Phi _{2}  \label{7.3}
\end{equation}
where $\varepsilon =-1$ for $\Bbb{R}^{4}$ and $\varepsilon =1$ for $\Bbb{R}%
^{2,2}$. Both the reductions (\ref{7.2}) and (\ref{7.3}) are admissible by
all equations of the DSII hierarchy if one chooses $A_{n}$ and $D_{n}$ as in
(\ref{A.13}).

In our case we have different linear problems for $\psi _{1}$, $\varphi _{1}$
and $\psi _{2}$, $\varphi _{2}$. To have the same equations for $p$ it is
enough to take 
\begin{equation}
\begin{array}{l}
A_{2n-1}=\partial _{\overline{z}}^{2n-1}+\ldots , \\ \vspace{-5mm} \\
A_{2n}=i\partial _{\overline{z}}^{2n}+\ldots ,
\end{array}
\begin{array}{l}
D_{2n-1}=\partial _{z}^{2n-1}+\ldots , \\ \vspace{-5mm} \\
D_{2n}=-i\partial _{z}^{2n}+\ldots
\end{array}
\label{7.4}
\end{equation}
in equations for $\psi _{1}$, $\varphi _{1}$ and 
\begin{equation}
\begin{array}{l}
A_{2n-1}=\partial _{\overline{z}}^{2n-1}+\ldots , \\ \vspace{-5mm} \\
A_{2n}=-i\partial _{\overline{z}}^{2n}+\ldots ,
\end{array}
\begin{array}{l}
D_{2n-1}=\partial _{z}^{2n-1}+\ldots , \\ \vspace{-5mm} \\
D_{2n}=i\partial _{z}^{2n}+\ldots
\end{array}
\label{7.5}
\end{equation}
in the case of $\psi _{2}$, $\varphi _{2}$. In particular, one gets 
\begin{equation}
\begin{array}{l}
p_{t_{2}}=i\left( p_{z\,z}+p_{\overline{z\,}\overline{z}}+up\right) \;\;, 
\\ \vspace{-5mm} \\ 
u_{z\overline{z}}=-2\varepsilon |p|_{z\,z}^{2}-2\varepsilon |p|_{\overline{%
z\,}\overline{z}}^{2}
\end{array}
\label{7.6}
\end{equation}
and 
\begin{equation}
\begin{array}{l}
\psi _{1t_{2}}=i\left( \partial _{\overline{z}}^{2}+w_{1}\right) \psi
_{1}+i\left( p_{z}-p\partial _{z}\right) \varphi _{1}\;\;, \\ 
\vspace{-5mm} \\
\varphi _{1t_{2}}=-i\varepsilon \left( \overline{p}_{\overline{z}}-\overline{%
p}\partial _{\overline{z}}\right) \psi _{1}-i\left( \partial
_{z}^{2}+w_{2}\right) \varphi _{1}
\end{array}
\label{7.7}
\end{equation}
while 
\begin{equation}
\begin{array}{l}
\psi _{2t_{2}}=-i\left( \partial _{\overline{z}}^{2}+w_{1}\right) \psi
_{2}-i\left( \overline{p}_{z}-\overline{p}\partial _{z}\right) \varphi
_{2}\;\;, \\ \vspace{-5mm} \\
\varphi _{2t_{2}}=i\varepsilon \left( p_{\overline{z}}-p\partial _{\overline{%
z}}\right) \psi _{2}+i\left( \partial _{z}^{2}+w_{2}\right) \varphi _{2}
\end{array}
\label{7.8}
\end{equation}
where $w_{1z}=-2\varepsilon |p|_{\overline{z}}^{2}$, $w_{2\overline{z}%
}=2\varepsilon |p|_{z}^{2}$.

The $t_{3}$ deformation is given now by equation (\ref{A.16}) 
and the deformations
of $\psi _{1}$, $\varphi _{1}$ and $\psi _{2}$, $\varphi _{2}$ are given by (%
\ref{7.1}), (A.8) with the reduction (\ref{7.2}) and (\ref{7.3}),
respectively.

Thus, in the cases of $\Bbb{R}^{4}$ and $\Bbb{R}^{2,2}$ deformations of
surfaces are generated by the proper DSII equation (\ref{7.6}) and the
corresponding hierarchy. Properties of solutions of the DSII equation (\ref
{7.6}) are essentially different for different signs of $\varepsilon $.
Consequently, the properties of deformations of surfaces in $\Bbb{R}^{4}$
and $\Bbb{R}^{2,2}$ will differ too.

In both cases $C_{1}=\int |p|^{2}dx\,dy$ is the integral of motion for the
whole hierarchy. Hence, the Willmore functionals for surfaces immersed into $%
\Bbb{R}^{4}$ and $\Bbb{R}^{2,2}$ are invariant under deformations generated
by the DSII hierarchy.

For the time-like surfaces the linear problem (\ref{5.2}), (\ref{5.11}),
obviously, give rise to the DSI hierarchy. In the case of the time-like
surfaces in $\Bbb{R}^{3,1}$ the reduction is $q=\overline{p}$, and,
consequently, the deformations of surfaces given by (\ref{5.1}), (\ref{5.2})
are generated by the DSI hierarchy under the reduction $q=\overline{p}$. In
particular, one has equations (\ref{A.14})-(\ref{A.16}) 
with $\varepsilon =1$. Again,
the Willmore functional $W=2\int |p|^{2}d\xi \,d\eta $ is invariant under
all these deformations.

Time-like surfaces in $\Bbb{R}^{2,2}$ are associated with the two different
linear problems (\ref{5.11}). To have the same evolution equations for the
pair $p$ and $q$, one has to choose the deformations of $\psi $, $\varphi $, 
$\widetilde{\psi }$, $\widetilde{\varphi }$ in the following form 
\begin{equation}
\left( 
\begin{array}{l}
\psi _{\alpha } \\ 
\varphi _{\alpha }
\end{array}
\right) _{t_{n}}=\left( 
\begin{array}{ll}
A_{n} & B_{n} \\ 
C_{n} & D_{n}
\end{array}
\right) \left( 
\begin{array}{l}
\psi _{\alpha } \\ 
\varphi _{\alpha }
\end{array}
\right) \qquad \qquad \alpha =1,2\quad ,  \label{7.9}
\end{equation}
\begin{equation}
\left( 
\begin{array}{l}
\widetilde{\psi }_{\alpha } \\ 
\widetilde{\varphi }_{\alpha }
\end{array}
\right) _{t_{n}}=\left( 
\begin{array}{ll}
\widetilde{A}_{n} & \widetilde{B}_{n} \\ 
\widetilde{C}_{n} & \widetilde{D}_{n}
\end{array}
\right) \left( 
\begin{array}{l}
\widetilde{\psi }_{\alpha } \\ 
\widetilde{\varphi }_{\alpha }
\end{array}
\right) \qquad \qquad \alpha =1,2\quad ,  \label{7.10}
\end{equation}
where $A_{n}$, $B_{n}$, $C_{n}$, $D_{n}$, $\widetilde{A}_{n}$, $\widetilde{B}%
_{n}$, $\widetilde{C}_{n}$, $\widetilde{D}_{n}$ are differential operators
with real-valued coefficients and 
\begin{eqnarray*}
\widetilde{A}_{n} &=&(-1)^{n-1}A_{n}[p\leftrightarrow q]\quad ,\quad 
\widetilde{B}_{n}=(-1)^{n-1}B_{n}[p\leftrightarrow q]\quad , \\
\widetilde{C}_{n} &=&(-1)^{n-1}C_{n}[p\leftrightarrow q]\quad ,\quad 
\widetilde{D}_{n}=(-1)^{n-1}D_{n}[p\leftrightarrow q]\quad .
\end{eqnarray*}
For example, the $t_{2}-$flow for $\psi _{\alpha }$, $\varphi _{\alpha }$ is
given by (A.2), (A.6) with real $\alpha _{2}$ while for $\widetilde{\psi }%
_{\alpha }$, $\widetilde{\varphi }_{\alpha }$ it is given by (\ref{7.10})
with the opposite signs of $\alpha _{2}$ and substitution $p\leftrightarrow
q $. In both cases one has equation (A.5). The Willmore functional $W=2\int
qp\,d\xi \,d\eta $ is clearly an invariant of all these deformations.

Thus, though the deformations for surfaces in $\Bbb{R}^{4}$, $\Bbb{R}^{3,1}$
and $\Bbb{R}^{2,2}$ are governed by different nonlinear integrable
equations, they have the following common property.

\begin{theorem}
The DSII hierarchy generates integrable deformations of space-like surfaces
immersed into $\Bbb{R}^{4}$, $\Bbb{R}^{3,1}$ and $\Bbb{R}^{2,2}$ via the
generalized Weierstrass representations. The DSI hierarchy generates
integrable deformations of time-like surfaces in 
$\Bbb{R}^{3,1}$ and $\Bbb{R}^{2,2}$. In
all cases the Willmore functionals $W$ for surfaces are invariant under the
corresponding deformations ($W_{t_{n}}=0$).
\end{theorem}

DS hierarchy of integrable equations is well studied [{49}]-[{51}]. This
provides us a broad class of deformations of surfaces in $\Bbb{R}^{4}$, $%
\Bbb{R}^{3,1}$ and $\Bbb{R}^{2,2}$ given explicitly. Moreover, since the
inverse spectral transform method allows us to linearize the initial-value
problem 
\[
p(z,\overline{z},t_{n}=0)\,,q(z,\overline{z},t_{n}=0)\rightarrow p(z,%
\overline{z},t_{n})\,,q(z,\overline{z},t_{n}) 
\]
for soliton equations of the DS hierarchy (see \emph{e.g.} [{49}]-[{51}]),
then the generalized Weierstrass formulae allows us to linearize the
initial-value problem for the deformation of surfaces $\vec{X}(z,\overline{z}%
,0)\rightarrow \vec{X}(z,\overline{z},t_{n})$. In virtue of all that, the
deformations generated by the DS hierarchy can be referred as integrable one.

Higher integrals of motion for the DS hierarchy are also certain 
functionals
on surfaces invariant under deformations generated by the DS hierarchy.
Since the Willmore functional $W$ is invariant under the conformal
transformations in four-dimensional spaces, then it is quite natural to
suggest

\begin{conjecture}
Higher integrals of motion for the DS hierarchy are functionals on
surfaces in $\Bbb{R}^{4}$, $\Bbb{R}^{3,1}$ and $\Bbb{R}^{2,2}$ which are
invariant under conformal transformations in these spaces.
\end{conjecture}

For tori in $\Bbb{R}^{3}$ an analogous conjecture has been proved in
[17].

Deformations of $\psi $, $\varphi $, $p$ allows us to find the deformation
equations for coordinates $\overrightarrow{X}$ and other geometrical
quantities. In the case of the space $\Bbb{R}^{4}$, using (\ref{3.15}), 
(\ref{7.6}), (\ref{7.7}), (\ref{7.8}), one obtains 
\begin{equation}
\begin{array}{l}
X_{1t_{2}}=
\frac{1}{2i}\left[ \overline{\psi }_{2}\overleftrightarrow{%
\partial }\overline{\psi }_{1}-\varphi _{2}\overleftrightarrow{\partial }%
\varphi _{1}-c.c.\right] \quad \quad , \\ \vspace{-5mm} \\
X_{2t_{2}}=\frac{1}{2}\left[ \overline{\psi }_{2}\overleftrightarrow{%
\partial }\overline{\psi }_{1}+\varphi _{2}\overleftrightarrow{\partial }%
\varphi _{1}+c.c.\right] \quad \quad , \\ \vspace{-5mm} \\
X_{3t_{2}}=\frac{1}{2i}\left[ \varphi _{2}\overleftrightarrow{\partial }%
\overline{\psi }_{1}+\overline{\psi }_{2}\overleftrightarrow{\partial }%
\varphi _{1}-c.c.\right] +I_{1}\quad , \\ \vspace{-5mm} \\
X_{4t_{2}}=\frac{1}{2}\left[ \varphi _{2}\overleftrightarrow{\partial }%
\overline{\psi }_{1}-\overline{\psi }_{2}\overleftrightarrow{\partial }%
\varphi _{1}+c.c.\right] +I_{2}\quad ,
\end{array}
\label{7.11}
\end{equation}
where $f\overleftrightarrow{\partial }g=f\,\partial g-g\,\partial f$ and 
\begin{eqnarray*}
I_{1} &=&\frac{1}{2i}\int_{\Gamma }\left[ \left( \overline{w}%
_{1}-w_{2}\right) \left( \overline{\psi }_{1}\varphi _{2}-\varphi _{1}%
\overline{\psi }_{2}\right) dz-c.c.\right] \quad , \\
I_{2} &=&\frac{1}{2}\int_{\Gamma }\left[ \left( \overline{w}%
_{1}-w_{2}\right) \left( \overline{\psi }_{1}\varphi _{2}+\varphi _{1}%
\overline{\psi }_{2}\right) dz+c.c.\right] \quad .
\end{eqnarray*}
The deformations (\ref{7.11}) of the coordinates can be decomposed into the
normal and tangential parts 
\begin{equation}
\overrightarrow{X}_{t_{2}}=a\overrightarrow{N}_{1}+b\overrightarrow{N}_{2}+c%
\overrightarrow{X}_{z}+\overline{c}\overrightarrow{X}_{\overline{z}}
\quad . \label{7.12}
\end{equation}
Using (\ref{3.3}), (\ref{3.7}), (\ref{3.8-3.9}) and (\ref{7.11}), one gets 
\begin{eqnarray*}
a &=&\frac{i}{2}\left( \left| \varphi _{1}\right| ^{2}\left| \varphi
_{2}\right| ^{2}u_{1}u_{2}\right) ^{-1/2}\times  \\
&&\times \left\{ \left[ u_{1}\left( \varphi _{1}\psi _{1}\overline{\varphi }%
_{2}\overline{\partial }\overline{\varphi }_{1}+\left| \varphi _{2}\right|
^{2}\varphi _{1}\partial \overline{\psi }_{2}\right) -c.c.\right]
-1\leftrightarrow 2\right\} +\Delta _{a}\quad ,\\
b &=&\frac{1}{2}\left( \left| \varphi _{1}\right| ^{2}\left| \varphi
_{2}\right| ^{2}u_{1}u_{2}\right) ^{-1/2}\times  \\
&&\times \left\{ \left[ u_{1}\left( \varphi _{2}\psi _{2}\overline{\varphi }%
_{1}\overline{\partial }\overline{\varphi }_{2}-\left| \varphi _{2}\right|
^{2}\varphi _{1}\partial \overline{\psi }_{2}\right) +c.c.\right]
+1\leftrightarrow 2\right\} +\Delta _{b}\quad ,\\
c & =& \frac{i}{2}\left( u_{1}u_{2}\right) ^{-1}\left\{
u_{1}u_{2z}-u_{2}u_{1z}\right\} +\Delta _{c}
\end{eqnarray*}
where
\begin{eqnarray*}
\Delta _{a} &=&\left( \left| \varphi _{1}\right| ^{2}\left| \varphi
_{2}\right| ^{2}u_{1}u_{2}\right) ^{-1/2}\times  \\
&&\times \left[ I_{1}\left( \func{Re}(\psi _{1}\varphi _{1}\overline{\psi }%
_{2}\overline{\varphi }_{2})-\left| \varphi _{1}\right| ^{2}\left| \varphi
_{2}\right| ^{2}\right) +I_{2}\func{Im}(\psi _{1}\varphi _{1}\overline{\psi }%
_{2}\overline{\varphi }_{2})\right] \quad ,\\
\Delta _{b} &=&\left( \left| \varphi _{1}\right| ^{2}\left| \varphi
_{2}\right| ^{2}u_{1}u_{2}\right) ^{-1/2}\times  \\
&&\times \left[ I_{1}\func{Im}(\psi _{1}\varphi _{1}\overline{\psi }_{2}%
\overline{\varphi }_{2})-I_{2}\left( \func{Re}(\psi _{1}\varphi _{1}%
\overline{\psi }_{2}\overline{\varphi }_{2})+\left| \varphi _{1}\right|
^{2}\left| \varphi _{2}\right| ^{2}\right) \right] \quad ,\\
\Delta _{c} &=& 
\frac{1}{2}\left( u_{1}u_{2}\right) ^{-1}\left[ I_{1}\left( \psi
_{1}\overline{\varphi }_{1}+\psi _{2}\overline{\varphi }_{2}\right)
-iI_{2}\left( \psi _{1}\overline{\varphi }_{1}-\psi _{2}\overline{\varphi }%
_{2}\right) \right] \quad .
\end{eqnarray*}

\section{Explicit deformations of surfaces}
\setcounter{equation}{0}

All explicit solutions of the DS hierarchy provide us deformations of
surfaces given by explicit formulae. We will present here two classes of
deformations for the space-like surfaces immersed in
$\Bbb{R}^{4}$ and $\Bbb{R}^{2,2}$
generated by the two basic classes of solutions of the DSII
hierarchy. We will give only final formulae omitting all calculations which
can be found in [49-54].

The first class is given by solutions which are parametrized by arbitrary
functions of a single variable. The corresponding $\psi _{\alpha }$, $%
\varphi _{\alpha }$and $p$ are given by [50,13] 

\begin{equation}
\begin{array}{ll}
\psi ^{(\alpha )} & =e^{-\lambda \overline{z}+\sum_{n=2}\lambda ^{n}t_{n}}+
\\ 
& +\frac{1}{2\pi i}\int \int_{\Bbb{C}}\frac{d\mu \wedge d\overline{\mu }}{%
\mu -\lambda }\sum_{k,m=1}^{N}\left[ \xi _{k}^{(\alpha )}\left( 1+A^{(\alpha
)}\right) _{km}^{-1}g_{m}^{(\alpha )}(\mu )\Gamma ^{-1}(\mu )\Gamma (\lambda
)\right] _{11}\,, \\ 
&  \\ 
\varphi ^{(\alpha )} & =\frac{1}{2\pi i}\int \int_{\Bbb{C}}\frac{d\mu \wedge
d\overline{\mu }}{\mu -\lambda }\sum_{k,m=1}^{N}\left[ \xi _{k}^{(\alpha
)}\left( 1+A^{(\alpha )}\right) _{km}^{-1}g_{m}^{(\alpha )}(\mu )\Gamma
^{-1}(\mu )\Gamma (\lambda )\right] _{21}\quad , \\ 
&  \\ 
p & =\frac{1}{2\pi i}\int \int_{\Bbb{C}}\frac{d\mu \wedge d\overline{\mu }}{%
\mu }\sum_{k,m=1}^{N}\left[ \xi _{k}^{(1)}\left( 1+A^{(1)}\right)
_{km}^{-1}g_{m}^{(1)}(\mu )\Gamma ^{-1}(\mu )\right] _{12}\quad ,
\end{array}
\label{8.1}
\end{equation}
where $\alpha =1,2$ , 
\begin{equation}
A_{km}(z,\overline{z})=\frac{1}{2\pi i}\int \int_{\Bbb{C}}d\lambda \wedge d%
\overline{\lambda }\int \int_{\Bbb{C}}\frac{d\mu \wedge d\overline{\mu }}{%
\mu -\lambda }g_{k}^{(\alpha )}(\mu )\Gamma ^{-1}(\mu )\Gamma (\lambda
)f_{m}^{(\alpha )}(\lambda )  \label{8.2}
\end{equation}
\begin{equation}
\begin{array}{l}
f_{k}^{(\alpha )}\left( \lambda ,\overline{\lambda }\right) =\left( 
\begin{array}{ll}
f_{1k}^{(\alpha )}(\lambda ,\overline{\lambda }) & \varepsilon \overline{%
f_{2k}^{(\alpha )}(-\lambda ,-\overline{\lambda })} \\ 
g_{2k}^{(\alpha )}(\lambda ,\overline{\lambda }) & \overline{f_{1k}^{(\alpha
)}(-\lambda ,-\overline{\lambda })}
\end{array}
\right) \quad , \\ 
\\ 
g_{k}^{(\alpha )}\left( \lambda ,\overline{\lambda }\right) =\left( 
\begin{array}{ll}
g_{1k}^{(\alpha )}(\lambda ,\overline{\lambda }) & \varepsilon \overline{%
f_{2k}^{(\alpha )}(-\overline{\lambda },-\lambda )} \\ 
g_{2k}^{(\alpha )}(\lambda ,\overline{\lambda }) & \overline{g_{1k}^{(\alpha
)}(-\overline{\lambda },-\lambda )} 
\end{array} 
\right) \quad ,
\end{array}
\label{8.3}
\end{equation}
\begin{eqnarray}
\xi _{m}^{(\alpha )}(z,\overline{z})&=&
\int \int_{\Bbb{C}}d\lambda \wedge d%
\overline{\lambda }\,\Gamma (\lambda )f_{m}^{(\alpha )}(\lambda )\quad , 
\nonumber \\ 
\\ 
\eta _{k}^{(\alpha )}(z,\overline{z})&=&\frac{1}{2\pi i}\int \int_{\Bbb{C}%
}d\mu \wedge d\overline{\mu }\,g_{k}^{(\alpha )}(\mu )\,\Gamma ^{-1}(\mu )
\label{8.4}
\end{eqnarray}
and 
\[
\Gamma \left( \lambda \right) =\left( 
\begin{array}{ll}
e^{-\lambda \overline{z}+\sum_{n=2}\lambda ^{n}t_{n}} & 0 \\ 
0 & e^{\lambda z-\sum_{n=2}\lambda ^{n}t_{n}}
\end{array}
\right) 
\]
and 
\[
f_{k}^{(2)}=\overline{f}_{k}^{(1)}\qquad ,\qquad g_{k}^{(2)}=\overline{g}%
_{k}^{(1)}\quad . 
\]
The $2\times 2$ matrices $\xi _{m}$ and $n_{k}$ are, in fact, of the form 
\begin{equation}
\xi _{m}^{(\alpha )}=\left( 
\begin{array}{ll}
\overline{\xi }_{m1}^{(\alpha )}(z) & \xi _{m2}^{(\alpha )}(z) \\ 
\varepsilon \overline{\xi }_{m2}^{(\alpha )}(z) & \xi _{m1}^{(\alpha )}(z)
\end{array}
\right) \quad ,\quad \eta _{m}^{(\alpha )}=\left( 
\begin{array}{ll}
\overline{\eta }_{k1}^{(\alpha )}(z) & \eta _{k2}^{(\alpha )}(z) \\ 
\varepsilon \overline{\eta }_{k2}^{(\alpha )}(z) & \eta _{k1}^{(\alpha )}(z)
\end{array}
\right)
\end{equation}
where $\xi _{m\beta }^{(\alpha )}(z)$, $\eta _{k\beta }^{(\alpha )}(z)$ ($%
m,k=1,\ldots $ , $\alpha ,\beta =1,2$ ) are arbitrary holomorphic functions.
They are, in essence, the complex Fourier transform of the functions $%
f_{\alpha k}(-\overline{\lambda })$ and $g_{\alpha k}(-\overline{\lambda })$%
. So the solution (\ref{8.1}) is parametrized by 4N arbitrary holomorphic
functions.

Consequently, the generalized Weierstrass formulae (\ref{3.15}) give us a
family of immersed surfaces parametrized by 4N arbitrary holomorphic
functions. Varying the times $t_{n}$, one gets integrable deformations of
these surfaces.

The multisoliton solutions form the second class of exact solutions. In this
case [52] 
\begin{equation}
\begin{array}{l}
p=-2\varepsilon \sum_{i=1}^{N}\overline{A}_{2i}(z)\,e^{-\overline{\lambda
_{i}\,z}} \\ 
\\ 
\overline{\psi }_{1}=\,e^{\lambda \,z}+\sum_{i=1}^{N}\frac{A_{1i}(z)}{%
\lambda _{i}-\lambda }\,e^{(\lambda -\lambda _{i})\,z} \\ 
\\ 
\overline{\varphi }_{1}=\sum_{i=1}^{N}\frac{\overline{A_{2i}(z)}}{\lambda
_{i}-\lambda }\,e^{(\lambda -\lambda _{i})\,z}
\end{array}
\label{8.7}
\end{equation}
where the column $\mathbf{A}_{i}=\left( 
\begin{array}{l}
A_{1i} \\ 
A_{2i}
\end{array}
\right) $ obey the system of equations 
\begin{equation}
\begin{array}{r}
\left( 
\begin{array}{l}
e^{\lambda _{n}z} \\ 
0
\end{array}
\right) -
\sum_{j=1, j\neq n} ^{N}\frac{\mathbf{A}_{j}}{\lambda
_{n}-\lambda _{j}}e^{(\lambda _{n}-\lambda _{i})z}=\left( z+\mu _{n}\right) 
\mathbf{A}_{n}+\nu _{n}\left( 
\begin{array}{ll}
0 & 1 \\ 
\varepsilon & 0
\end{array}
\right) \overline{\mathbf{A}}_{n}\quad \quad , \\ 
n=1,...,N\quad \quad ,
\end{array}
\label{8.8}
\end{equation}
where $\lambda _{i}$, $\mu _{n}$, $\nu _{n}$ are arbitrary complex
constants. For $\psi _{2}$ and $\varphi _{2}$ one has similar expressions
with the substitution $\lambda _{i}\rightarrow \widetilde{\lambda }_{i}$, $%
\mu _{n}\rightarrow \widetilde{\mu }_{n}$, $\nu _{n}\rightarrow \widetilde{%
\nu }_{n}$ which should be chosen so that 
\begin{equation}
\sum_{i=1}^{n}\overline{\widetilde{A}}_{2i}\,\exp \left( -\overline{%
\widetilde{\lambda }_{i}z}\right) =\sum_{i=1}^{n} {A}_{2i}\,\exp
\left( -\lambda _{i}z \right) \quad .  \label{8.9}
\end{equation}
This is the condition that from the linear problem for $\psi _{2}$, $\varphi
_{2}$ one gets $\overline{p}$ instead of $p$ at the problem for $\psi _{1}$, 
$\varphi _{1}$. It can be shown [52] that 
\begin{equation}
|p|^{2}=-4\varepsilon \left[ \log \det D\right] _{z\,\overline{z}}
\label{8.10}
\end{equation}
where the $2N\times 2N$ matrix $D$ is given by 
\begin{equation}
D=\left( 
\begin{array}{ll}
M & Q \\ 
\overline{Q} & \overline{M}
\end{array}
\right)  \label{8.11}
\end{equation}
and 
\begin{eqnarray}
&&M_{jj}=z+\mu _{j}+2i\,\lambda _{j}\,t_{2}\quad \quad ,  \nonumber \\
&&M_{ij}=(\lambda _{i}-\lambda _{j})^{-1}e^{(\lambda _{i}-\lambda
_{j})\,z}\quad \quad \quad i\neq j\quad ,  \label{8.12} \\
&&Q=diag\left[ \nu _{i}\,e^{-it_{2}\left( \lambda _{j}^{2}+\overline{\lambda 
}_{j}^{2}\right) }\right] \quad \quad \quad j=1,...,N  \nonumber
\end{eqnarray}
Here we are restricted to the solutions of equation (\ref{7.6}). Finally one
has [52] 
\begin{equation}
\int |p|^{2}dx\,dy=-4\pi \varepsilon \,N  \quad .\label{8.13}
\end{equation}
The simplest solution is ($N=1$)
\begin{equation}
\begin{array}{l}
p=2\overline{\nu }\frac{\exp \left[ \lambda _{1}z-\overline{\lambda }_{1}%
\overline{z}+i(\lambda _{1}^{2}+\overline{\lambda }_{1}^{2})\,t_{2}\right] }{%
|z+2i\lambda _{1}t_{2}+\mu |^{2}-\varepsilon |v|^{2}} \quad , \\ 
\\ 
\overline{\psi }=e^{\lambda _{1}z}-\left( \lambda -\lambda _{1}\right)
^{-1}e^{\lambda \,z}\left[ \log \left( |z+2i\lambda _{1}t_{2}+\mu
|^{2}-\varepsilon |v|^{2}\right) \right] _{z} \quad , \\ 
\\ 
\overline{\varphi }=\frac{\varepsilon }{2}\left( \lambda -\lambda
_{1}\right) ^{-1}e^{\lambda \,z\,}p \quad .
\end{array}
\label{8.14}
\end{equation}
These solutions give rise to surfaces and their deformations of soliton
type. For them the value of the Willmore functional is 
\begin{equation}
W=16\pi N\quad .  \label{8.15}
\end{equation}
Note that it does not depend on the sign of $\varepsilon $.

For the DSI equation there is another interesting class of explicit
solutions. For these solutions, called dromions, the functions $p$ and $q$
decrease exponentially fast in all directions on the plane $(\xi ,\eta )$
(see \emph{e.g.} [53,49,51]). For general dromion solution of equation
(\ref{A.14}) one has [54] 
\begin{equation}
|p|^{2}=-4\left[ \log \det \left( 1-\varepsilon A\right) \right] _{\xi \eta }
\label{8.16}
\end{equation}
where the rectangular matrix A is of the form 
\begin{equation}
A=\varrho \beta \varrho ^{\dagger }\overline{\alpha }  \label{8.17}
\end{equation}
and 
\begin{equation}
\begin{array}{l}
\left( \rho \right) _{ij}=\rho _{ij}\quad \quad ,\quad \quad (\rho ^{\dagger
})_{ij}=\overline{\rho }_{ji}\quad , \\ 
\\ 
\beta _{ij}=\int_{-\infty }^{\eta }Y_{i}(\eta ^{\prime },t_{2})\,Y_{j}(\eta
^{\prime },t_{2})\,d\eta ^{\prime }\quad , \\ 
\\ 
\alpha _{ij}=\int_{-\infty }^{\xi }X_{i}(\xi ^{\prime },t_{2})\,X_{j}(\xi
^{\prime },t_{2})\,d\xi ^{\prime }\quad .
\end{array}
\label{8.18}
\end{equation}
Here $X_{i}$ and $Y_{i}$ are arbitrary solutions of equations 
\begin{equation}
\begin{array}{l}
iX_{i\,t_{2}}+X_{i\xi \xi }+u_{2}(\xi ,t_{2})X_{i}=0\quad , \\ 
\\ 
iY_{i\,t_{2}}+Y_{i\eta \eta }+u_{1}(\eta ,t_{2})Y_{i}=0
\end{array}
\label{8.19}
\end{equation}
and $u_{1}$ and $u_{2}$ are arbitrary functions. Choosing $u_{1}$ and $u_{2}$
as the reflectionless potentials in the Schroedinger equations (\ref{8.19}),
one gets multidromion solutions. For these solutions of the DSI equation one
has [54] 
\begin{equation}
\int |p|^{2}d\xi \,d\eta =-\varepsilon \log \det \left( 1-\varepsilon
\varrho \varrho ^{\dagger }\right)  \quad . \label{8.30}
\end{equation}
As a result the corresponding time-like surfaces in the Minkowski space $%
\Bbb{R}^{3,1}$ (the formulae (\ref{5.1}),(\ref{5.2})) have the following
value of the Willmore functional ($\varepsilon =1$) 
\begin{equation}
W=-2\log \det \left( 1-\varrho \varrho ^{\dagger }\right)  
\quad . \label{8.21}
\end{equation}
This class of surfaces could be of interest for the theory of classical
strings in the Minkowski space.

The formulae (\ref{8.15}) and (\ref{8.21}) 
demonstrate us that the method of the
inverse scattering transform provides us the technique to calculate the
value of the Willmore functional for rather complicated surfaces in 4D
spaces.

\section{Particular classes of surfaces and their deformations}
\setcounter{equation}{0}

Here we consider some special classes of surfaces which have simple
geometric characterization.

We start with minimal surfaces. In the four-dimensional spaces they are
characterized by the condition $\overrightarrow{H}=0$. Using the expressions
for the mean curvature vector found in sections 3-5, we conclude that in all
cases $\overrightarrow{H}=0$ if the potentials $p$ and $q$ in DS-linear
problems vanish. So one has

\begin{corollary}
The formulae of section 3-5 with $p=q=0$ give us the Weierstrass
representations for minimal surfaces in 4D spaces. Space-like surfaces are
parametrized by four arbitrary holomorphic functions (two holomorphic and
two anti-holomorphic). Formulae for time-like surfaces contain two arbitrary
functions of one variable ($\xi $) and two arbitrary functions of another
variable ($\eta $).
\end{corollary}

From the formulae for $\vec{H}^{2}$ follows also that surfaces with zero
length of the mean curvature vector ($\vec{H}^{2}=0$) coincide with the
minimal surfaces in the cases of space-like surfaces in $\Bbb{R}^{4}$, $\Bbb{%
R}^{2,2}$ and time-like surfaces in $\Bbb{R}^{3,1}$. In contrast, this
condition is less restrictive in the cases of space-like surfaces in $\Bbb{R}%
^{3,1}$ and time-like surfaces in $\Bbb{R}^{2,2}$: it is sufficient that $%
p=0 $ (or $q=0$).

Superminimal immersions form a subclass of the minimal ones for which in
addition the condition $\vec{X}_{zz}\cdot \vec{X}_{zz}=0$ is satisfied 
(see \textit{e.g.} [29]). For
surfaces in $\Bbb{R}^{4}$ and $\Bbb{R}^{2,2}$, using (\ref{3.15}) and (\ref
{4.1}), one gets 
\begin{equation}
\vec{X}_{zz}\cdot \vec{X}_{zz}=\varepsilon \left( \overline{\psi }%
_{1z}\varphi _{1}-\overline{\psi }_{1}\varphi _{1z}\right) \left( \overline{%
\psi }_{2z}\varphi _{2}-\overline{\psi }_{2}\varphi _{2z}\right) =0
\label{9.1}
\end{equation}
where $\varepsilon =-1$ for $\Bbb{R}^{4}$ and $\varepsilon =1$ for $\Bbb{R}%
^{2,2}$ .

\begin{corollary}
The immersions of surface into $\Bbb{R}^{4}$ and $\Bbb{R}^{2,2}$ 
given by the
Weierstrass representations (\ref{3.15})- (\ref{3.16}) and (\ref{4.1})-(\ref
{4.2}), respectively, are superminimal if $p=0$ and $\varphi _{1}=a_{1}%
\overline{\psi }_{1}$ (or $\varphi _{2}=a_{2}\overline{\psi }_{2}$) where $%
a_{1}$ (or $a_{2}$) is an arbitrary constant.
\end{corollary}

Analogous results hold for the other cases.

The Weierstrass representations for superminimal immersions could be useful
for an analysis of the problems discussed in [29].

Next geometrically interesting class of surfaces correspond to the constant
mean curvature (\emph{i.e.} $\vec{H}^2=const.$). In the case of $\Bbb{R}^{4}$
one, obviously, has

\begin{corollary}
Surfaces of the constant length of the mean curvature vector ($\vec{H}%
^{2}=const$) conformally immersed into $\Bbb{R}^{4}$ are generated by the
formulae (\ref{3.15}) where $\psi _{\alpha }$, $\varphi _{\alpha }$ ($\alpha
=1,2$) obey to the system of equations 
\begin{eqnarray}
\psi _{\alpha z}&=&\frac{1}{2}\exp (i\theta_{\alpha})
\,\sqrt{\vec{H}^{2}\left( |\psi
_{1}|^{2}+|\varphi _{1}|^{2}\right) \left( |\psi _{2}|^{2}+|\varphi
_{2}|^{2}\right) }\;\varphi _{\alpha }\quad , \nonumber \\ 
\varphi _{\alpha \overline{z}}&=&
-\frac{1}{2}\exp (-i\theta_{\alpha})\,
\sqrt{\vec{H}%
^{2}\left( |\psi _{1}|^{2}+|\varphi _{1}|^{2}\right) \left( |\psi
_{2}|^{2}+|\varphi _{2}|^{2}\right) }\;\psi _{\alpha }
\quad , \quad \alpha =1,2  \nonumber \\
&&
\label{9.2}
\end{eqnarray}
where $\theta_1=-\theta_2=\theta$ and $\theta (z,\overline{z)}$ 
is an arbitrary function.
\end{corollary}

In the space-like case in $\Bbb{R}^{2,2}$ the surfaces of the constant $\vec{%
H}^{2}$ are generated by the formula (\ref{4.1}) where $\psi _{\alpha }$, $%
\varphi _{\alpha }$ obey (\ref{9.2}) with obvious changes of signs. Similar
situation take place for time-like surfaces in $\Bbb{R}^{3,1}$ while for
space-like surfaces in $\Bbb{R}^{3,1}$ and time-like surfaces in $\Bbb{R}%
^{2,2}$ the functions 
$\psi _{\alpha }$, $\varphi _{\alpha }$ obey the system (\ref{4.8})
and (\ref{5.11}) with the constraints 
\begin{equation}
qp=-\frac{1}{4}\vec{H}^{2}\,\left| \psi _{1}\varphi _{2}-\psi _{2}\varphi
_{1}\right| ^{2}  \label{9.3}
\end{equation}
and 
\begin{equation}
qp=-\frac{1}{4}\vec{H}^{2}\,\left[ \left( \psi _{1}\varphi _{2}-\psi
_{2}\varphi _{1}\right) \left( \widetilde{\psi }_{1}\widetilde{\varphi }_{2}-%
\widetilde{\psi }_{2}\widetilde{\varphi }_{1}\right) \right] ^{-1}
\label{9.4}
\end{equation}
respectively.

Other special classes of surfaces in 4D spaces are associated with the
reductions of the corresponding linear problems. The constraint $\overline{p}%
=p$ in the linear systems (\ref{3.16}), (\ref{4.2}), (\ref{5.2}) gives rise
to special classes of surfaces in $\Bbb{R}^{4}$, $\Bbb{R}^{2,2}$, and $\Bbb{%
R}^{3,1}$. These cases for spaces $\Bbb{R}^{4}$ and $\Bbb{R}^{2,2}$ have
been considered earlier in [35]. The corresponding 
integrable deformations are generated by the modified Veselov-Novikov 
(VN) hierarchy.

Analogously, the constraint $p=q$ for the Weierstrass representations (\ref
{4.7})- (\ref{4.8}) and (\ref{5.10})- (\ref{5.11}) gives rise to special
classes of space-like surfaces in $\Bbb{R}^{3,1}$ and time-like surfaces in $%
\Bbb{R}^{2,2}$, respectively. The integrable deformations are generated
correspondingly by the real modified Nizhnik-Veselov-Novikov (NVN) 
hierarchy (reduction of the
DSII and DSI hierarchies under the constraint $p=q$).

Another interesting class of surfaces given by the formulae (\ref{4.7})-(\ref
{4.8}) and (\ref{5.10})-(\ref{5.11}) is associated with the reduction $p=1$.
In this case the basic linear system is equivalent to the equation 
\begin{equation}
\psi _{\xi \eta }=q\psi  \label{9.5}
\end{equation}
($\xi =z$, $\eta =\overline{z}$ for $\Bbb{R}^{3,1}$ space). Since $\psi
_{\xi }=\varphi $ the Weierstrass formulae are simplified. For instance, the
formulae (\ref{4.7})-(\ref{4.8}) become 
\begin{equation}
\begin{array}{l}
X^{1}+iX^{2}=\int_{\Gamma }(\overline{\psi }_{2}\psi _{1z^{\prime
}}\,dz^{\prime }+\psi _{1}\overline{\psi }_{2\overline{z}^{\prime }}\,d%
\overline{z}^{\prime })\quad , \\ \vspace{-5mm} \\
X^{1}-iX^{2}=\int_{\Gamma }(\overline{\psi }_{1}\psi _{2z^{\prime
}}\,dz^{\prime }+\psi _{2}\overline{\psi }_{1\overline{z}^{\prime }}\,d%
\overline{z}^{\prime })\quad , \\ \vspace{-5mm} \\
X^{3}+X^{4}=\int_{\Gamma }(\overline{\psi }_{1}\psi _{1z^{\prime
}}\,dz^{\prime }+\psi _{1}\overline{\psi }_{1\overline{z}^{\prime }}\,d%
\overline{z}^{\prime })\quad , \\ \vspace{-5mm} \\
X^{3}-X^{4}=-\int_{\Gamma }(\overline{\psi }_{2}\psi _{2z^{\prime
}}\,dz^{\prime }+\psi _{2}\overline{\psi }_{2\overline{z}^{\prime }}\,d%
\overline{z}^{\prime })
\end{array}
\label{9.6}
\end{equation}
where 
\begin{equation}
\psi _{\alpha z\,\overline{z}}=q\psi _{\alpha }\qquad ,\qquad \alpha =1,2
\label{9.7}
\end{equation}
$q$ is a real-valued function and 
\begin{equation}
ds^{2}=\left| w\right| ^{2}dz\,d\overline{z}  \label{9.8}
\end{equation}
where $w$ is the Wronskian 
\begin{equation}
w=\left| 
\begin{array}{ll}
\psi _{1} & \psi _{1z} \\ 
\psi _{2} & \psi _{2z}
\end{array}
\right| \quad .  \label{9.9}
\end{equation}
Integrable deformations of this class of surfaces in $\Bbb{R}^{3,1}$ are
generated by the VN equation (\ref{A.17}) and the whole VN hierarchy.

For time-like surfaces in $\Bbb{R}^{2,2}$ the corresponding Weierstrass
representation is given by 
\begin{equation}
\begin{array}{l}
X^{1}=\frac{1}{2}\int_{\Gamma }[(\psi _{1\xi ^{\prime }}\psi _{4}+\psi
_{2\xi ^{\prime }}\psi _{3})\,d\xi ^{\prime }+(\psi _{1}\psi _{4\eta
^{\prime }}+\psi _{2}\psi _{3\eta ^{\prime }})\,d\eta ^{\prime }]\quad , 
\\ \vspace{-5mm} \\ 
X^{2}=\frac{1}{2}\int_{\Gamma }[(\psi _{1\xi ^{\prime }}\psi _{3}-\psi
_{2\xi ^{\prime }}\psi _{4})\,d\xi ^{\prime }+(\psi _{1}\psi _{3\eta
^{\prime }}-\psi _{2}\psi _{4\eta ^{\prime }})\,d\eta ^{\prime }]\quad , 
\\ \vspace{-5mm} \\ 
X^{3}=\frac{1}{2}\int_{\Gamma }[(\psi _{2\xi ^{\prime }}\psi _{3}-\psi
_{1\xi ^{\prime }}\psi _{4})\,d\xi ^{\prime }+(\psi _{2}\psi _{3\eta
^{\prime }}-\psi _{1}\psi _{4\eta ^{\prime }})\,d\eta ^{\prime }]\quad , 
\\ \vspace{-5mm} \\ 
X^{4}=\frac{1}{2}\int_{\Gamma }[(\psi _{1\xi ^{\prime }}\psi _{3}+\psi
_{2\xi ^{\prime }}\psi _{4})\,d\xi ^{\prime }+(\psi _{1}\psi _{3\eta
^{\prime }}+\psi _{2}\psi _{4\eta ^{\prime }})\,d\eta ^{\prime }]
\end{array}
\label{9.10}
\end{equation}
where the real-valued functions $\psi _{1}$, $\psi _{2}$, $\psi _{3}$, $\psi
_{4}$ satisfy the equation 
\begin{equation}
\psi _{\alpha \eta \xi }=q\psi _{\alpha }\qquad ,\qquad \alpha =1,2,3,4
\label{9.11}
\end{equation}
with real-valued $q$. Here we denote $\widetilde{\varphi }_{1}$, $\widetilde{%
\varphi }_{2}$ from (\ref{5.11}) as $\widetilde{\varphi }_{1}=\psi _{3}$, $%
\widetilde{\varphi }_{2}=\psi _{4}$. The induced metric is 
\begin{equation}
ds^{2}=(\psi _{1\xi }\psi _{2}-\psi _{1}\psi _{2\xi })\,(\psi _{3\eta }\psi
_{4}-\psi _{3}\psi _{4\eta })\,d\xi \,d\eta  \quad . \label{9.12}
\end{equation}
Integrable deformations are generated by the Nizhnik equation (equation
(\ref{A.17}) with real-valued $\xi $ and $\eta $) and by the whole Nizhnik
hierarchy.

The Weierstrass formulae (\ref{9.6}) and (\ref{9.10}) represent surfaces in
the spaces $\Bbb{R}^{3,1}$ and $\Bbb{R}^{2,2}$ via four solutions of the
two-dimensional Schroedinger (Moutard, perturbed string) equation (\ref{9.11}%
) (or (\ref{9.7})). So, they are the sort of 
four-dimensional extensions of the
Lelieuvre formula (\ref{2.12})-(\ref{2.13}). But now surfaces generated are
parametrized by the minimal lines instead of asymptotic lines in (\ref{2.12}%
)-(\ref{2.13}). The fact that integrable deformations of surfaces in 4D
spaces represented by (\ref{9.6}) and (\ref{9.10}) and of affine surfaces
generated by the Lelieuvre formula (\ref{2.12})-(\ref{2.13}) are governed by
the same NVN hierarchy suggests a certain connection between these classes
of surfaces.

\section{One dimensional reductions}
\setcounter{equation}{0}

All the Weierstrass formulae presented above have one natural special case
when the potentials (coefficients) in the basic linear systems depend only
on one variable, say $x$ ($x=\func{Re}z$ or $x=\frac{\xi +\eta }{2}$ ). Let
us consider space-like surfaces conformally immersed in $\Bbb{R}^{4}$, $\Bbb{%
R}^{2,2}$ or $\Bbb{R}^{3,1}$. Let $p=p(x)$, $q=q(x)$. In this case solutions 
$\psi $, $\varphi $ of the basic linear system are of the form 
\begin{equation}
\psi =\chi (x)\exp (\lambda y)\qquad ,\qquad \varphi =\widetilde{\chi }%
(x)\exp (\lambda y)  \label{10.1}
\end{equation}
where $\lambda $ is an arbitrary complex parameter and $y=\func{Im}(z)$ or $%
y=(\xi -\eta )/2$. For the second linear problem the functions $\psi _{2}$, $%
\varphi _{2}$ also have the form (\ref{10.1}) but, in general, with another
parameter $\mu $. Correspondingly, the basic linear system is reduced to the
following one-dimensional one 
\begin{equation}
\left( 
\begin{array}{l}
\chi \\ 
\widetilde{\chi }
\end{array}
\right) _{x}=i\lambda \left( 
\begin{array}{cc}
1 & 0 \\ 
0 & -1
\end{array}
\right) \left( 
\begin{array}{l}
\chi \\ 
\widetilde{\chi }
\end{array}
\right) +2\left( 
\begin{array}{cc}
0 & p \\ 
q & 0
\end{array}
\right) \left( 
\begin{array}{l}
\chi \\ 
\widetilde{\chi }
\end{array}
\right) \quad .  \label{10.2}
\end{equation}
The AKNS hierarchy of 1+1-dimensional integrable equations (\ref{A.19}) is
associated with the spectral problem (\ref{10.2}). The properties of this
system and, consequently, of the one-dimensional reduction of the
Weierstrass representation are different for the cases of real $\lambda $
and imaginary $\lambda $.

Let us consider first the case of real $\lambda $ and $\mu $. For the real $%
\lambda $ and $\mu $ the formulae (\ref{3.15}), (\ref{4.1}) and (\ref{4.7})
imply that 
\begin{equation}
\overrightarrow{X}(x,y)=\overrightarrow{Y}(x)\,\exp [(\lambda +\mu )y]\quad .
\label{10.3a}
\end{equation}
Equivalently 
\begin{equation}
\left( \partial _{\overline{z}}-\partial _{z}\right) \overrightarrow{X}%
=i(\lambda +\mu )\overrightarrow{X}\quad .  \label{10.3}
\end{equation}
Hence 
\begin{equation}
\left( \partial _{\overline{z}}-\partial _{z}\right) \overrightarrow{X}\cdot
\left( \partial _{\overline{z}}-\partial _{z}\right) \overrightarrow{X}%
=-(\lambda +\mu )^{2}\overrightarrow{X}\cdot \overrightarrow{X}\quad .
\label{10.3b}
\end{equation}
Since for conformal immersion 
\begin{equation}
\partial _{z}\overrightarrow{X}\cdot \partial _{z}\overrightarrow{X}%
=\partial _{\overline{z}}\overrightarrow{X}\cdot \partial _{\overline{z}}%
\overrightarrow{X}=0\quad ,  \label{10.4}
\end{equation}
equation (\ref{10.3}) implies 
\begin{equation}
\partial _{z}\overrightarrow{X}\cdot \partial _{\overline{z}}\overrightarrow{%
X}=\frac{1}{2}(\lambda +\mu )^{2}\overrightarrow{X}\cdot \overrightarrow{X}%
\quad .  \label{10.5}
\end{equation}
Due to (\ref{3.4}), (\ref{3.6}), (\ref{4.3}) and (\ref{4.9}) 
\begin{equation}
\partial _{z}\overrightarrow{X}\cdot \partial _{\overline{z}}\overrightarrow{%
X}=\frac{1}{2}\det \Phi _{1}\det \Phi _{2}  \label{10.6}
\end{equation}
where $\det \Phi _{1}$ and $\det \Phi _{2}$ are the determinants of the $%
2\times 2$ matrix solutions of the form (\ref{4.13}), (\ref{4.14}), (\ref
{4.15}) for the systems (\ref{3.16}), (\ref{4.1}), (\ref{4.8}). Using (\ref
{10.1}) and (\ref{10.2}), one gets from (\ref{10.5}) and (\ref{10.6}) the
relation 
\begin{equation}
\det \widehat{\Phi }_{1}\det \widehat{\Phi }_{2}=(\lambda +\mu )^{2}%
\overrightarrow{Y}(x)\cdot \overrightarrow{Y}(x)  \label{10.7}
\end{equation}
where $\widehat{\Phi }_{1}$ and $\widehat{\Phi }_{2}$ are $2\times 2$ matrix
solutions of the system (\ref{10.2}) and of the corresponding second system
, respectively. It is easy to check that for the problem (\ref{10.2}) $(\det 
\widehat{\Phi })_{x}=0$. Without lost of generality one can put $\det 
\widehat{\Phi }_{1}=\det \widehat{\Phi }_{2}=1$ (that is typical for the
one-dimensional spectral problems). Consequently, the relation (\ref{10.7})
gives 
\begin{equation}
\overrightarrow{Y}(x)\cdot \overrightarrow{Y}(x)=\frac{1}{(\lambda +\mu )^{2}%
}  \quad . \label{10.8}
\end{equation}
Surface with coordinates $\overrightarrow{X}$ of the form (\ref{10.2}) with
real $\lambda +\mu $ is of a cone type. It is obtained by homothety of the
curve with coordinate $\overrightarrow{Y}(x)$ by the factor $\exp [(\lambda
+\mu )y]$. The relation (\ref{10.8}) shows that this curve lies on the
sphere $\Bbb{S}^{3}$ of the radius $\frac{1}{(\lambda +\mu )}$ at the
original space $\Bbb{R}^{4}$or on the hyperboloids for spaces $\Bbb{R}^{2,2}$
and $\Bbb{R}^{3,1}$.

Integrable deformations of surfaces and, consequently, of the curves with
the coordinates $\overrightarrow{Y}(x)$ are given by the AKNS hierarchy. In
particular, in the case of the space $\Bbb{R}^{4}$ ($q=\overline{p}$) one
has an integrable motion of curves on the sphere $\Bbb{S}^{3}$ of the radius 
$\frac{1}{(\lambda +\mu )}$ which is governed by the NLS hierarchy. So we 
reproduce the result of the paper [55] about integrable motion of curves on 
$\Bbb{S}^{3}$. The case $\lambda +\mu =0$ corresponds to the integrable
motion of space curves in $\Bbb{R}^{3}$ [56].

For pseudo-Euclidean spaces $\Bbb{R}^{2,2}$ and $\Bbb{R}^{3,1}$ one has
integrable motions of curves on the three-dimensional hyperboloids governed
by the NLS hierarchy ($\varepsilon =1$) and by the real AKNS hierarchy ($%
\Bbb{R}^{3,1}$).

Thus the one-dimensional limit of the Weierstrass formulae and corresponding
deformations reproduces the results for integrable motion of curves in
three-dimensional spaces.

A different situation arises in the case of pure imaginary $\lambda =i\omega 
$ in (\ref{10.1}). Indeed, for instance, for surfaces immersed into $\Bbb{R}%
^{4}$ with such $\psi _{\alpha }$ and $\varphi _{\alpha }$, one has (see eq.
(\ref{3.3})) 
\begin{equation}
\begin{array}{ll}
\partial \left( X^{1}+iX^{2}\right) =-\widetilde{\chi }_{1}\widetilde{\chi }%
_{2}e^{i\omega _{+}y}\quad , & \partial \left( X^{1}-iX^{2}\right) =%
\overline{\chi }_{1}\overline{\chi }_{2}e^{-i\omega _{+}y}\quad , \\ 
\vspace{-5mm} \\
\partial \left( X^{3}+iX^{4}\right) =\widetilde{\chi }_{1}\overline{\chi }%
_{2}e^{i\omega _{-}y}\quad , & \partial \left( X^{3}-iX^{4}\right) =%
\overline{\chi }_{1}\widetilde{\chi }_{2}e^{-i\omega _{-}y}
\end{array}
\label{10.9}
\end{equation}
where $\omega _{\pm }=\omega _{1}\pm \omega _{2}$ and $\chi _{i},\widetilde{%
\chi }_{i}$ depend only on $x$. The induced metric is 
\begin{equation}
ds^{2}=A(x)\,dz\,d\overline{z}=\left( |\chi _{1}|^{2}+|\widetilde{\chi }%
_{1}|^{2}\right) \left( |\chi _{2}|^{2}+|\widetilde{\chi }_{2}|^{2}\right)
\,dz\,d\overline{z}\quad .  \label{10.10}
\end{equation}
So this reduction gives rise to a surface of ''revolution'' in $\Bbb{R}^{4}$%
. The linear problem (\ref{10.2}) in this case and the corresponding NLS
hierarchy are studied in great details (see \emph{e.g.} [49,57,58]). We will
apply all these results to surfaces of revolution in $\Bbb{R}^{4}$ in a
separate paper.

\section*{Appendix. DS hierarchy.}
\setcounter{equation}{0}
\renewcommand{\thesection}{\Alph{section}}
\addtocounter{section}{-9}

Here we present some basic known facts about the DS equation and the DS
hierarchy. They are associated (see \emph{e.g.} [49-51]) with the following
two-dimensional linear system (Dirac equation) 
\begin{equation}
\begin{array}{l}
\psi _{\xi }=p\varphi \\ 
\varphi _{\eta }=q\psi
\end{array}
 \label{A.1}
\end{equation}
where $\psi $, $\varphi $, $p$, $q$ are, in general, complex-valued
functions of the independent variables $\xi $, $\eta $ which can be either
real or complex. In soliton theory this system is known as the
Davey-Stewartson (DS) linear problem. An infinite hierarchy of nonlinear
differential equations associated with (A.1) is referred as the DS
hierarchy. It arises as the compatibility condition of (A.1) with the
systems [{49}-{51}] 
\begin{equation}
\begin{array}{l}
\psi _{t_{n}}=A_{n}\psi +B_{n}\varphi \;\;, \\ 
\varphi _{t_{n}}=C_{n}\psi +D_{n}\varphi 
\end{array}
  \label{A.2}
\end{equation}
where times $t_{n}$ are new (deformation) variables and $A_{n}$, $B_{n}$, $%
C_{n}$, $D_{n}$ are differential operators of $n-$th order. At $n=1$ one
gets the linear system 
\begin{equation}
\begin{array}{l}
p_{t_{1}}=\alpha p_{\eta }+\gamma p_{\xi }\;\;, \\ 
q_{t_{1}}=\gamma q_{\xi }+\alpha q_{\eta }
\end{array}
\label{A.3}
\end{equation}
where $\alpha $, $\gamma $ are arbitrary constants. The corresponding
operators in (\ref{A.3}) are 
\begin{equation}
A_{1}=\alpha \partial _{\eta }\;,\;B_{1}=\gamma p\;,\;C_{1}=\alpha
q\;,\;D_{1}=\gamma \partial _{\xi }  
\label{A.4}
\end{equation}
Higher equations are nonlinear ones. 
At $n=2$ one has the system [{49}]-[{51}] 
\begin{equation}
\begin{array}{l}
p_{t_{2}}=\alpha _{2}\left( p_{\xi \xi }+p_{\eta \eta }+up\right) ;\;, \\ 
q_{t_{2}}=-\alpha _{2}\left( q_{\xi \xi }+q_{\eta \eta }+uq\right) ;\;, \\ 
u_{z\overline{z}}=-2(pq)_{\xi \xi }-2(pq)_{\eta \eta }
\end{array}
\label{A.5}
\end{equation}
where $\alpha _{2}$ is an arbitrary constant. For the system (A.5) 
\begin{equation}
\begin{array}{l}
A_{2}=\alpha _{2}\left( \partial _{\eta }^{2}+w_{1}\right) \;, \\ 
C_{2}=-\alpha _{2} ( q_{\eta }-q\partial _{\eta } ) \;,
\end{array}
\begin{array}{l}
B_{2}=\alpha _{2}\left( p_{\xi }-p\partial _{\xi }\right) \;, \\ 
D_{2}=-\alpha _{2} ( \partial _{\xi }^{2}+w_{2} )
\end{array}
\label{A.6}
\end{equation}
where 
\[
w_{1\xi }=-2(pq)_{\eta }\;\;,\;\;w_{2\eta }=-2(pq)_{\xi
}\;\;,\;\;u=w_{1}+w_{2}\;\;. 
\]
For the $t_{3}$ deformations one has (see \emph{e.g.} [50]) 
\begin{equation}
\begin{array}{l}
p_{t_{3}}=p_{\xi \xi \xi }+p_{\eta \eta \eta }+3p_{\xi }\partial _{\eta
}^{-1}(pq)_{\xi }+3p_{\eta }\partial _{\xi }^{-1}(pq)_{\eta }+3p\partial
_{\eta }^{-1}(qp_{\xi })_{\xi }+3q\partial _{\xi }^{-1}(qp_{\eta })_{\eta
}\;, \\ 
\vspace{-7mm}
\\
q_{t_{3}}=q_{\xi \xi \xi }+q_{\eta \eta \eta }+3q_{\xi }\partial _{\eta
}^{-1}(pq)_{\xi }+3q_{\eta }\partial _{\xi }^{-1}(pq)_{\eta }+3q\partial
_{\eta }^{-1}(pq_{\xi })_{\xi }+3q\partial _{\xi }^{-1}(pq_{\eta })_{\eta
}\;\;.
\end{array}
\label{A.7}
\end{equation}
In this case 
\begin{equation}
\begin{array}{l}
A_{3}=\partial _{\eta}^{3}+3\left[ \partial _{\xi
}^{-1}(pq)_{\eta }\right] \partial _{\eta }+3\left[ \partial _{\xi
}^{-1}(qp_{\eta })_{\eta }\right] \;, \\ \vspace{-5mm}\\
B_{3}=-p\partial _{\xi }^{2}+p_{\xi }\partial _{\xi }-p_{\xi \xi }-3p\left[
\partial _{\eta }^{-1}(pq)_{\xi }\right] \;, \\ \vspace{-5mm} \\
C_{3}=-q\partial _{\eta }^{2}+q_{\eta }\partial _{\eta }-q_{\eta \eta
}-3q\left[ \partial _{\xi }^{-1}(pq)_{\eta }\right] \;, \\ \vspace{-5mm} \\
D_{3}=\partial _{\xi }^{3}+3\left[ \partial _{\eta }^{-1}(pq)_{\xi }\right]
\partial _{\xi }+3\left[ \partial _{\eta }^{-1}(pq_{\xi })_{\xi }\right] \;.
\end{array}
\label{A.8}
\end{equation}
For $t_{n}-$evolution the operators $A_{n},D_{n}$ are of the order $n$ and $%
B_{n},C_{n}$ are of the order $n-1$.

The system (A.5) is referred as the DS system. The DS hierarchy for real $%
\xi $ and $\eta $ is referred as the DSI hierarchy while at the case $\xi =z$%
, $\eta =\overline{z}$ it is called DSII hierarchy.

The DS hierarchy, known also as the two-component KP hierarchy, can be
written in different 
compact forms. First, within the Sato approach it is equivalent
to the infinite system of equations [59] 
\begin{equation}
\frac{\partial Q}{\partial t_{n\alpha }}=Q\left( Q^{-1}H_{\alpha }\partial
^{n}Q\right)_{-} \quad ,\quad \alpha=1,2 \quad , \quad n=1,2,\ldots 
\label{A.9}
\end{equation}
where the pseudo-differential operator $Q$ is $Q=1+w_{1}\partial
^{-1}+w_{2}\partial ^{-2}+\ldots $, $w_{k}$ are $2\times 2$ matrices, $%
\partial =\partial _{t_{1}}$, the matrices $H_{\alpha }$ form a basis of the
diagonal $2\times 2$ matrices and $\pounds _{-}=\pounds -\pounds _{+}$ where 
$\pounds _{+}$ is the differential part of the operator $\pounds $.

Second, the DS hierarchy can be written with the use of the bilocal
recursion operator $L(x,y,t_{n)}$ as [60] ($\xi =x+y$ , $\eta =x-y$ or $%
\xi =x+iy$ , $\eta =x-iy$ ) 
\begin{equation}
\frac{\partial P}{\partial t_{n_{\alpha}}}=\Delta \,\beta _{\alpha }L^{n}\left(
H_{\alpha }P^{\prime }-PH_{\alpha }\right) \qquad ,\qquad 
\alpha=1,2 \quad, \quad n=1,2,... 
\label{A.10}
\end{equation}
where 
\[
P^{\prime }=P(x,y^{\prime },t_{n})=\left( 
\begin{array}{cc}
0 & p(x,y^{\prime },t) \\ 
q(x,y^{\prime },t) & 0
\end{array}
\right) \qquad , 
\]
$\beta _{\alpha }$ are arbitrary constants, the bilocal recursion operator $%
L $ acts as follows 
\[
L\chi (x,y^{\prime },y^{\prime })=(\partial _{x}+\sigma _{3}\partial _{y})%
\widetilde{\chi }+\partial _{y^{\prime }}\widetilde{\chi }^{\prime }\sigma
_{3}-P\,d^{-1}(P\widetilde{\chi }-\widetilde{\chi }P^{\prime })_{D}+d^{-1}(P%
\widetilde{\chi }-\widetilde{\chi }P^{\prime })_{D}P^{\prime }(x,y)\, 
\]
where $[\sigma _{3},\widetilde{\chi }]=\chi $, $d=\partial _{x}+\sigma
_{3}(\partial _{y}+\partial _{y^{\prime }})$ and $Z_{D}$ means the diagonal
part of the $2\times 2$ matrix $Z$. The bilocality of the recursion operator
is the principal feature of the 2+1 integrable systems.

The equations from the DS hierarchy are integrable by the inverse spectral
transform method (see [49-51]). They possess all remarkable properties
typical for the 2+1-dimensional soliton equations, namely: there are 
infinite classes of solutions given by explicit formulae 
(solutions with functional parameters,
multisoliton solutions, periodic solutions), infinite symmetry algebra,
Darboux and Backlund transformations, Hamiltonian and Lagrangian structures
etc. . They have an infinite set of integrals of motion $C_{n}$ (\emph{i.e. }%
$\frac{\partial C_{n}}{\partial t_{m}}=0$, $n,m=1,2,\ldots $). The simplest
of them is 
\begin{equation}
C_{1}=\int qp\,d\xi \,d\eta 
\label{A.12}
\end{equation}
while higher integrals of motion are non-local. Emphasize that (\ref{A.12}) 
is the
integral of motion for the whole DS hierarchy ($\frac{\partial C_{1}}{%
\partial t_{m}}=0$, $m=1,2,\ldots $).

The DS hierarchy contains different sub-hierarchies associated with
different specializations of $p$ and $q$ (reductions). The most important
reduction is $q=\varepsilon \overline{p}$ where $\varepsilon =\pm 1$. This
reduction is compatible with all equations of the DS hierarchy if one
chooses 
\begin{equation}
\begin{array}{ll}
A_{2n-1}=\partial _{\eta }^{2n-1}+\ldots \quad , & D_{2n-1}=\partial _{\xi
}^{2n-1}+\ldots \quad , \\ 
A_{2n}=\pm i\partial _{\eta }^{2n}+\ldots \quad , & D_{2n}=\mp i\partial
_{\xi }^{2n}+\ldots \qquad \qquad (n=1,2,3,4).
\end{array}
\label{A.13}
\end{equation}
In particular, the $t_{2}-$flow takes a form 
\begin{equation}
\begin{array}{l}
p_{t_{2}}=i\left( p_{\xi \xi }+p_{\eta \eta }+up\right) \;\;, \\ 
u_{\xi \eta }=-2\varepsilon |p|_{\xi \xi }^{2}-2\varepsilon |p|_{\eta \eta
}^{2}
\end{array}
\label{A.14}
\end{equation}
for which 
\begin{equation}
\begin{array}{l}
\psi _{t_{2}}=i\left( \partial _{\eta }^{2}+w_{1}\right) \psi +i\left(
p_{\xi }-p\partial _{\xi }\right) \varphi \quad , \\ 
\varphi _{t_{2}}=-i\varepsilon \left( \overline{p}_{\eta }-\overline{p}%
\partial _{\eta }\right) \psi -i\left( \partial _{\xi
}^{2}+w_{2}\right) \varphi
\end{array}
\label{A.15}
\end{equation}
where $w_{1\xi }=-2\varepsilon |p|_{\eta }^{2}$, $w_{2\eta }=2\varepsilon
|p|_{\xi }^{2}$.

The $t_{3}$ deformation is given now by the equation 
\begin{equation}
p_{t_{3}}=p_{\xi \xi \xi }+p_{\eta \eta \eta }+3\varepsilon p_{\xi }\partial
_{\eta }^{-1}(|p|_{\xi }^{2})+3\varepsilon p_{\eta }\partial _{\xi
}^{-1}(|p|_{\eta }^{2})+3\varepsilon p\partial _{\eta }^{-1}(\overline{p}%
p_{\xi })_{\xi }+3\varepsilon p\partial _{\xi }^{-1}(\overline{p}p_{\eta
})_{\eta } 
\label{A.16}
\end{equation}
where the evolution of $\psi $ and $\varphi $ is given by (A.2), (A.8) with 
$q=\varepsilon \overline{p}$.

Equation (\ref{A.14}) 
is the proper DS equation which appears in hydrodynamics and
plasma physics\ (see [49]). At $\varepsilon =-1$ and $\varepsilon =1$ one
has the defocusing and focusing cases, respectively.

Deeper reductions of the DS hierarchy are associated with constraints $p=\pm
q$ or $p=1$. Both of them are admissible only by equations with odd times .
Under the reduction $p=1$ the system (A.7) is converted into the equation 
\begin{equation}
q_{t_{3}}=q_{\xi \xi \xi }+q_{\eta \eta \eta }+3\left[ q\partial _{\eta
}^{-1}(q_{\xi })\right] _{\xi }+3\left[ q\partial _{\xi }^{-1}(q_{\eta
})\right] _{\eta }  
\label{A.17}
\end{equation}
and the linear system (A.1) is equivalent to 
\begin{equation}
\psi _{\xi \eta }=q\psi \quad .  
\label{A.18}
\end{equation}
Equation (\ref{A.17}) 
is the well-known Nizhik-Veselov-Novikov (NVN) equation. 
Under this reduction the DS hierarchy is converted into the NVN hierarchy.

In the one-dimensional limit $p_{\xi }=p_{\eta }$, $q_{\xi }=q_{\eta }$, 
\emph{i.e. }$p=p(x)=p\left( \frac{\xi +\eta }{2}\right) $ , $q=q(x)=q\left( 
\frac{\xi +\eta }{2}\right) $, the DS hierarchy is reduced to the AKNS
hierarchy of the 1+1-dimensional integrable system. Using the recursion
operator one can represent this hierarchy in the form (see \emph{e.g.} [49]) 
\begin{equation}
\left( 
\begin{array}{c}
q \\ 
-p
\end{array}
\right) _{t_{n}}=\alpha _{n}L^{n}\left( 
\begin{array}{l}
q \\ p
\end{array}
\right) \qquad ,\qquad n=1,2,\ldots  
\label{A.19}
\end{equation}
where $\alpha _{n}$ are constants and the recursion operator $L$ is 
\begin{equation}
L=\left( 
\begin{array}{ll}
\partial _{x}-2q\partial _{x}^{-1}p\qquad & 2q\partial _{x}^{-1}q \\ 
-2p\partial _{x}^{-1}p & -\partial _{x}+2p\partial _{x}^{-1}q
\end{array}
\right) \quad .  
\label{A.20}
\end{equation}
The simplest equations from the AKNS hierarchy look like 
\begin{equation}
\begin{array}{l}
q_{t_{2}}=\alpha _{2}\left( q_{xx}-2q^{2}p\right) \quad , \\ 
p_{t_{2}}=-\alpha _{2}\left( p_{xx}-2p^{2}q\right)
\end{array}
 \label{A.21}
\end{equation}
and 
\begin{equation}
\begin{array}{l}
q_{t_{3}}=\alpha _{3}\left[ q_{xxx} -6 p q q_x \right] \quad , \\ 
p_{t_{3}}=\alpha _{3}\left[ p_{xxx}-6 q p p_x \right] \quad .
\end{array}
 \label{A.22}
\end{equation}
Under the reduction $q=\varepsilon \overline{p}$ the system (\ref{A.21}) 
becomes ($\alpha _{2}=i$) 
\begin{equation}
ip_{t_{2}}=p_{xx}-2\varepsilon |p|^{2}p \label{A.23}
\end{equation}
that is the famous nonlinear Schroedinger (NLS) equation. At the case $p=-1$
($\alpha _{3}=-1$) the system (\ref{A.22}) is reduced to the celebrated Korteweg
de Vries (KdV) equation 
\begin{equation}
q_{t_{3}}=-q_{xxx}-6qq_{x}  \label{A.24}
\end{equation}
while under the reduction $q=p$ ($\alpha _{3}=-1$) one gets the modified KdV
equation 
\begin{equation}
p_{t_{3}}=p_{xxx}-6p^{2}p_{x}  \quad . \label{A.25}
\end{equation}
Equations (\ref{A.14}), (\ref{A.17}) and (\ref{A.16}) 
($\overline{p}=p$) are the
2+1-dimensional integrable generalizations of the NLS, KdV and mKdV
equations, respectively.

\textbf{Acknowledgements.} One of the authors (B.G.K.)
appreciates very much the discussions with E. Ferapontov, F.
Pedit, U. Pinkall and T. Willmore. B.G.K. is also grateful to D. Fairlie 
for attracting his attention to the paper [46].

\hfill 
\hfill
\hfill

%\newpage

\begin{centerline}
\textbf{REFERENCES}
\end{centerline}

\begin{enumerate}
\item  \label{r1} G. Darboux, \textit{Lecons sur la th\'{e}orie des surfaces
et les applications geometriques du calcul infinitesimal}, \emph{t. 1-4},
Gauthier-Villars, Paris, 1877-1896.

\item  \label{r2} L. Bianchi, \textit{Lezioni di Geometria Differenziale}, 
\textit{2nd ed.}, Spoerri, Pisa, 1902.

\item  \label{r3} L.P. Eisenhart, \textit{A treatise on the differential
geometry of Curves and Surfaces}, Dover, New York, 1909.

\item  T.J. Willmore, \textit{Total curvature in Riemannian Geometry}, Ellis
Horwood, New York, 1982.

\item  S.S. Chern, in: \textit{Differential Geometry and Complex Analysis}
(Eds., I. Chavel, H.M. Farkas), Springer, Berlin, 1985.

\item  \label{r4} S.T. Yau (Ed.), \textit{Seminar on Differential Geometry},
Princeton Univ., Princeton, 1992.\label{r5}

\item  \label{r6} D. Nelson, T. Piran and S. Weinberg (Eds.), \textit{%
Statistical mechanics of membranes and surfaces}, Worls Scientific,
Singapore, 1989.

\item  \label{r7} D.J. Gross, T. Piran and S. Weinberg (Eds.), \textit{Two
dimensional quantum gravity and random surfaces}, World Scientific,
Singapore, 1992.

\item  \label{r8} F. David, P. Ginsparg and Y. Zinn-Justin (Eds.), \textit{%
Fluctuating geometries in Statistical Mechanics and Field Theory}, Elsevier
Science, Amsterdam, 1996.

\item  R.C. Brower, D.A. Kessler, J. Koplik and H. Levine, \textit{%
Geometrical models of interface evolution}, Phys. Rev. A 29: 1335-1342,
(1984).

\item  K. Nakayama and M. Wadati, \textit{The motion of surfaces}, J. Phys.
Soc. Japan 62:1895-1901, (1993).

\item  \label{r10} B.G. Konopelchenko, \textit{Multidimensional integrable
systems and dynamics of surfaces in space}, preprint of Institute of
Mathematics, Taipei, August 1993; in \textit{National Workshop on Nonlinear
Dynamics}, (M. Costato, A. Degasperis and M. Milani, Eds.), Ital. Phys.
Society, Bologna, Italy, 1995, pp. 33-40.

\item  B.G. Konopelchenko, \textit{Induced surfaces and their integrable
dynamics}, Stud. Appl. Math. 96:9-51, (1996).

\item  R. Carroll and B.G. Konopelchenko, \textit{Generalized
Weierstrass-Enneper inducing, conformal immersion and gravity}, Int. J.
Modern Physics A 11:1183-1216, (1996).

\item  \label{r14} B.G. Konopelchenko and I. Taimanov, \textit{Generalized
Weierstrass formulae, soliton equations and Willmore surfaces}, preprint N.
187, Univ. Bochum, (1995).

\item  \label{r15} B.G. Konopelchenko and I. Taimanov, \textit{Constant mean
curvature surfaces via an integrable dynamical system}, J. Phys. A: Math.
Gen. 29:1261-1265, (1996).

\item  \label{r16} I. Taimanov, \textit{Modified Novikov-Veselov equation
and differential geometry of surfaces}, Trans. Amer. Math. Soc., Ser. 2
,179:133-159, (1997).

\item  \label{r17} P.G. Grinevich and M.V. Schmidt, \textit{Conformal
invariant functionals of immersion of tori into }$\Bbb{R}$\textit{$^{3}$},
Journal of Geometry and Physics (to appear); preprint SFB288 N 291,
TU-Berlin, 1997.

\item  \label{r18} J. Richter, \textit{Conformal maps of a Riemann surface
into space of quaternions}, PH.D Thesis, TU-Berlin, 1997.

\item  \label{r19} I. Taimanov, \textit{Global Weierstrass representation
and its spectrum}, Uspechi Mat. Nauk. 52, N 6, 187-188, (1997).

\item  \label{r20} I. Taimanov, \textit{The Weierstrass representation of
spheres in }$\Bbb{R}$\textit{$^{3}$, the Willmore numbers and soliton spheres%
}, preprint SFB 288, N 302, TU-Berlin, 1998.

\item  I. Taimanov, \textit{Surfaces of revolution in terms of solitons,
Annals Global Anal. Geometry }15:419-435, (1997).

\item  B.G. Konopelchenko, \textit{On solutions of shape equations for
membranes and strings}, Phys. Lett. B 414:58-64, (1997).

\item  S. Matsutani, \textit{Constant mean curvature surfaces and the Dirac
operator}, J. Phys. A: Math. Gen. 30:4019-4024, (1997).

\item  S. Matsutani, \textit{Immersion anomaly of Dirac operator for
surfaces in }$\Bbb{R}$\textit{$^{3}$}, preprint physics/9707003 , (1997).

\item  \label{r22} S. Matsutani, \textit{On density state of quantized
Willmore surfaces-a way to quantized extrinsic string in $\Bbb{R}^{3}$}, 
J. Phys. A: Math. Gen. 31:3595-3606, (1998).

\item  \label{r23} B.G. Konopelchenko and G. Landolfi, \textit{On classical
string configurations}, Mod. Phys. Lett. A 12:3161-3168, (1997).

\item  B.G. Konopelchenko and U. Pinkall, \textit{Integrable deformations of
affine surfaces via Nizhnik-Veselov-Novikov equation}, Phys. Lett A
245:239-245, (1998).

\item  R.L. Bryant, \textit{Conformal and minimal immersions of compact
surfaces into the $4-$sphere}, J. Diff. Geom. 17:455-473, (1982).

\item  \label{r25} T. Friedrich, \textit{On surfaces in four-spaces}, Ann.
Glob. Analysis and Geometry 2:257-287, (1984).

\item  \label{r26} D.A. Hoffman and R. Osserman, \textit{The Gauss map of
surfaces in }$\Bbb{R}$\textit{$^{3}$ and }$\Bbb{R}$\textit{$^{4}$}, Proc.
London Math. Soc. (3) 50:27-56, (1985).

\item  \label{r27} D. Ferus, F. Pedit, U. Pinkall and I. Sterling, \textit{%
Minimal tori in }$\Bbb{S}$\textit{$^{4}$}, J. Reine Angew. Math. 429:1-47,
(1992).

\item  B.G. Konopelchenko, \textit{Weierstrass representations for surfaces
in 4D spaces and their integrable deformations via the DS hierarchy}, 
math.DG/9807129, (1998).

\item  P. Baird and J.C. Wood, \textit{Weierstrass representations for
harmonic morphisms on Euclidean spaces and surfaces}, dg-ga/9512010, (1995).

\item  \label{r29} B.G. Konopelchenko and G. Landolfi, \textit{Generalized
Weierstrass representation for surfaces in multidimensional Riemann spaces},
math.DG/ 9804144, (1998); J. Geometry and Physics (to appear).

\item  K. Kenmotsu, \textit{Weierstrass formula for surfaces of prescribed
mean curvature}, Math. Ann. 245:89-99, (1979).

\item  U. Abresh, \textit{Spinor representation of CMC surfaces}, Lecture at
Luminy, 1989.

\item  A.I. Bobenko, \textit{Surfaces in terms of }$2\times 2$\textit{\
matrices. Old and new integrable cases}, in: \textit{Harmonic maps and
integrable systems}, (Fordy A. and Wood J., Eds.), pp. 83-127, Vieweg, 1994.

\item  R. Kusner and N. Schmitt, \textit{The spinor representation of
surfaces in space}, dg-ga/9610005, (1996).

\item  T. Friedrich, \textit{On the spinor representation of surfaces in
Euclidean 3-spaces}, preprint SFB 288, N. 295, TU-Berlin, (1997).

\item  S. Matsutani, \textit{Dirac operators of a conformal surface immersed
in }$\Bbb{R}$\textit{$^{4}$: further generalized Weierstrass relation},
solv-int/9801006, (1998).

\item  F. Pedit and U. Pinkall, in preparation.

\item  G. Kamberov, F. Pedit and U. Pinkall, \textit{Bonnet pairs and
isothermic surfaces}, Duke Math. J. 92:637-644, (1998).

\item  X-F Cao and C. Tian, \textit{Integrable system and spacelike surfaces
with prescribed mean curvature in Minkowski 3-space}, preprint (1998).

\item  K. Akutagawa and S. Nishikawa, \textit{The Gauss map and spacelike
surfaces with prescribed mean curvature in Minkowski 3-space}, Tohoku Math.
J. 42:67-82, (1990).

\item  D.B. Fairlie and C.A. Manogue, \textit{Lorentz invariance and the
composite strings}, Phys. Rev. D 34:1832-1834, (1986).

\item  A.I. Bobenko, \textit{All constant mean curvature tori in }$\Bbb{R}%
^{3}$\textit{, }$\Bbb{S}^{3}$\textit{, }$\Bbb{H}^{3}$\textit{\ in terms of
theta-functions}, Math. Ann. 290:209-245, (1991).

\item  L.P. Eisenhart, \textit{Riemannian Geometry}, Princeton, 1926.

\item  M.J. Ablowitz and P.A. Clarkson, \textit{Solitons, nonlinear
evolution equations and inverse scattering}, Cambridge Univ. Press,
Cambridge, 1991.

\item  \label{r37} B.G. Konopelchenko, \textit{Introduction to
multidimensional integrable equations}, Plenum Press, New York, 1992.

\item  \label{r38} B.G. Konopelchenko, \textit{Solitons in multidimension},
World Scientific, Singapore, 1993.

\item  V.A. Arkadiev, A.K. Pogrebkov and M.C. Polivanov, \textit{Inverse
scattering transform method and soliton solutions for Davey-Stewartson II
equation}, Physica D 36:189-197, (1989).

\item  M. Boiti, L. Martina and F. Pempinelli, \textit{Multidimensional
localized solitons}, Chaos, Solitons and Fractals 5:2377-2417, (1995).

\item  P.M. Santini, \textit{Energy exchange of interacting coherent
structures in multidimensions}, Physica D 41:26-54, (1990).

\item  A. Doliwa and P.M. Santini, \textit{An elementary geometric
characterization of the integrable motions of a curve}, Phys. Lett A.
185:373-389, (1994).

\item  H. Hasimoto, \textit{A soliton on a vortex filament}, J. Fluid. Mech.
51:477-485, (1972).

\item  V.E. Zakharov, S.V. Manakov, S.P. Novikov and L.P. Pitaevsky, \textit{%
Theory of solitons}, Consultans Bureau, New York, 1984.

\item  L.D. Faddeev and L.A. Takhtajan, \textit{Hamiltonian methods in the
theory of solitons}, Springer-Verlag, Berlin, 1987.

\item  M. Jimbo and T. Miwa, \textit{Solitons and infinite-dimensional Lie
algebras}, Publ.\ RIMS, Kyoto Univ. 19:943-1001, (1983).

\item  A.S. Fokas and P.M. Santini, \textit{Bi-hamiltonian structures in
multidimension, I}, Commun. Math. Phys. 115:375-419, (1988); \textit{II},
Commun. Math. Phys. 116:449-479, (1988).
\end{enumerate}

\end{document}